\documentclass[12pt]{amsart}
\usepackage{amsmath}
\usepackage{amsfonts}
\usepackage{amssymb}
\usepackage{amsthm}
\usepackage{graphicx}
\usepackage{color}
\usepackage{eucal}
\usepackage{tikz}

\newtheorem{theorem}{Theorem}[section]

\newtheorem{corollary}[theorem]{Corollary}

\newtheorem{definition}[theorem]{Definition}
\newtheorem{example}[theorem]{Example}

\newtheorem{proposition}[theorem]{Proposition}
\newtheorem{remark}[theorem]{Remark}

\makeatletter
\def\ball{{B}}
\def\tto{\rightrightarrows}
\def\bx{\bar x}

\def\bz{\bar z}

\def\dom{\mathop{\rm dom}\nolimits}
\def\rge{\mathop{\rm rge}\nolimits}
\def\gph{\mathop{\rm gph}\nolimits}

\def\dist{\mathop{\rm dist}\nolimits}

\def\diam{\mathop{\rm diam}\nolimits}

\begin{document}

\subjclass{ 47J05  47J07  49J52 
49J53  58C15}


\title{A Quest for Simple and Unified Proofs in Regularity Theory: Perturbation Stability}



\author{Radek Cibulka} 
\address{NTIS - New Technologies for the Information Society and Department of Mathematics, Faculty of Applied Sciences, University of West Bohemia, Univerzitn\'i 8, 301 00 Plze\v{n}, Czech Republic}
\email{cibi@kma.zcu.cz}
 
\author{Tom\'a\v{s} Roubal}
\address{NTIS - New Technologies for the Information Society and Department of Mathematics, Faculty of Applied Sciences, University of West Bohemia, Univerzitn\'i 8, 301 00 Plze\v{n}, Czech Republic}
\email{roubalt@kma.zcu.cz}


%

\thanks{Both the authors were supported by the grant of GA\v CR no. 20-11164L}

%


\begin{abstract}
Ioffe's criterion and various reformulations of it  have become a~standard tool in proving theorems guaranteeing  metric regularity of a (set-valued) mapping. First, we demonstrate  that one should always use directly the so-called general criterion which follows, for example, from Ekeland's variational principle, and that there is no need to make a detour through the slope-based consequences of this general statement. Second, we argue that when proving perturbation stability results, in the spirit of Lyusternik-Graves theorem, there is no need to employ the concept of a lower semi-continuous envelope even in the case of an incomplete target space. The gist is to use the ``correct'' function to which  Ekeland's variational principle is applied; namely, the distance function to the graph of the set-valued mapping under consideration. This approach originates in the notion of graphical regularity introduced by L. Thibault, which is equivalent to the property of metric regularity. Our criteria cover  also both metric subregularity and metric semiregularity, which are weaker properties obtained  by fixing one of the points in the definition of metric regularity.
\end{abstract}

\keywords{metric regularity, graphical regularity, Ioffe criterion, perturbation stability}

\maketitle

\section{Introduction} \label{scIntroduction} Metric regularity, linear openness, and Aubin property of the inverse of a~given set-valued mapping are three \emph{equivalent} properties playing a fundamental role in modern variational analysis and have been broadly covered in several monographs \cite{AubinFrankowskaBook, BorweinZhuBook, DontchevBook, IoffeBook, MordukhovichBook, PenotBook}. We formulate all the results in terms of the first mentioned property and call it simply regularity. Although our notation is fairly standard, in case of any difficulty, the reader can consult the end of this introduction. We distinguish two main settings in which regularity statements for a set-valued mapping $F: X \rightrightarrows Y$ from a metric space $(X,d)$ into another metric space $(Y, \varrho)$ are proved:
\begin{itemize}
 \item[(S$_1$)] the graph of $F$ is (locally) complete;    
 \item[(S$_2$)] the domain space $X$ is complete and the graph of $F$ is (locally) closed.
\end{itemize} 
For example, $(S_1)$ occurs when both the domain space and the target space are complete and the graph of $F$ is closed in $X \times Y$, which often can be assumed, without any loss of generality, by using \cite[Proposition 2.40]{IoffeBook} provided that all the assumptions translate to the completions of the spaces and the graph.
There are several ways of proving sufficient conditions for regularity in the literature. Let us mention the most popular ones (from our point of view at least):

\smallskip
\paragraph*{\bf (A) Iterative processes.}
This (constructive) method finds a solution of the inclusion [equation] as the limit of a~sequence of iterates generated by a~process  resembling Newton's method or the contraction mapping iteration, see, for example, \cite{G}, \cite{D1996}, \cite[Proof I of Theorem 1A.1, p. 12]{DontchevBook}. One can say that these proofs are the longest ones but work  without any difficulty in the setting $(S_2)$, that is,  when $Y$ is not necessarily complete  \cite[Proof I of Theorem 5E.1, p. 309]{DontchevBook}.

\smallskip
\paragraph*{\bf (B) A direct application of Ekeland's variational principle (EVP)}
This approach is used frequently in \cite[Section 3.4]{AubinFrankowskaBook}, for example. We recall a weak form this principle. 
\begin{theorem}\label{Ekeland}
	Let $(X,d)$ be a complete metric space and an $x\in X$ be given.  Consider a lower semi-continuous  function $\varphi :X \to [0,\infty]$ such that $\varphi(x)<\infty$. Then there exists a point $u\in X$	such that $\varphi(u) +d(u,x )\leq \varphi(x)$	and
	\begin{equation} \label{eqEVP_01}
		\varphi(u) \leq \varphi(\hat{u})+ d(\hat{u},u)\quad \text{whenever}\quad \hat{u} \in X.
	\end{equation}
\end{theorem}
\noindent A  full version of this principle ensures that the inequality in \eqref{eqEVP_01} is strict provided that $\hat{u} \neq u$. However, we show that this is not needed for applications in regularity theory.  Theorem~\ref{Ekeland} is usually proved via an iterative process generating a (countable) sequence converging to a desired point. 
\smallskip
\paragraph*{\bf (C) Application of slope-based conditions.}  These conditions are usually derived either from the EVP directly or via a detour through error bounds \cite{HT08}. The concept of the (local strong) \emph{slope} of a~function $\varphi: X \to \mathbb{R}\cup\{\infty\}$ at a point  $x \in X$, with $\varphi(x) < \infty$, was introduced in De Giorgi, Marino, and Tosques \cite{DMT};  and this quantity is defined by 
$$
 | \nabla \varphi |(x) := 
 \limsup\limits_{u \to x, u \neq x} \frac{[\varphi(x) - \varphi(u)]^+}{d(u,x)}, 
$$
where $[a]^+ := \max\{a,0\}$, $a \in \mathbb{R}$; when $\varphi(x)=\infty$ then $| \nabla \varphi |(x) := \infty$. This approach, going back to works by D. Az\'e and J.N. Corvellec \cite{AC06},  can be found, for example, in \cite{ACN15},   \cite{HNT2014}, \cite{HT15}, \cite{HT08}; see \cite{IoffeBook} for more details concerning the historical background. Especially in the setting $(S_2)$, additional concepts such as lower semi-continuous envelopes are employed in the literature, which is tangling the proofs.

\smallskip
\paragraph*{\bf (D) Application of a fixed point theorem.}
A particular choice of the principle depends on the regularity statement we want to prove. In case of Graves' theorem \cite{G}, the usual Banach contraction mapping principle \cite[Theorem 1A.3]{DontchevBook} does the job, cf. 
\cite[Proof II of Theorem 1A.1, p. 19]{DontchevBook}; when considering the mappings $F + g$, with a \emph{smooth} single-valued part $g$ and a set-valued part $F$, we can choose  the contraction mapping principle for set-valued mappings \cite[Theorem 5E.2]{DontchevBook} proved by A.L. Dontchev and W.W. Hager in 1994. As the proofs of both the preceding fixed point theorems rely on the iterative procedure one can expect to obtain shorter proofs of the regularity statement in question. Similarly to {\bf (A)}, this approach works without any difficulty in the setting $(S_2)$,   see \cite[Proof of Theorem 5E.5]{DontchevBook}. Under certain compactness assumptions, Brouwer's or Kakutani's fixed point theorem can be used instead of contraction mapping principles. In finite-dimensions, one can obtain results that do not seem  to follow from the infinite-dimensional ones, see \cite{CR2022} for more details.  However, when considering the mappings $F + g$, with a \emph{non-smooth} single-valued part $g$ in infinite dimensions, it seems to be necessary to \emph{combine} Brouwer's fixed point theorem (or Glicksberg's extension of Kakutani's fixed point theorem) and the iterative procedure {\bf (A)} (or Ioffe's regularity criterion {\bf (E)} described below), see  \cite{CD2016, CDV2016, cf15}.
\smallskip
\paragraph*{\bf (E) Application of general Ioffe-type regularity criteria.} 
This approach, extensively used in \cite{IoffeBook}, is presented in this note. It encapsulates both the iterative procedure as in {\bf (A)} and Ekeland's variational principle from {\bf (B)}. In fact, one can use any of these two methods to get these criteria, see \cite[Proposition 2.1]{cf15} going back to \cite{FP1987} for the first and \cite[Theorem 2.45]{IoffeBook} going back to \cite{i2000} for the latter. Similarly to {\bf (B)} and {\bf (C)}, the key staring point is to choose an appropriate underlying lower semi-continuous function, see \cite[pp. 60--65]{IoffeBook} for several criteria formulated in terms of various functions described below. 
Given a single-valued continuous mapping $f$ from a complete metric space $(X,d)$ into a metric space $(Y,\varrho)$,   one applies Theorem~\ref{Ekeland} to the function $\varphi:= \varrho(y, f(\cdot))$, where $y$ is the right-hand side of the  equation $f(x) = y$ the solution of which we are searching for; hence no completeness assumption is imposed on the target space, cf. the setting $(S_2)$. There are several ways of handling set-valued mappings. Assume that a right-hand side $y \in Y$ is given and we want to find an $x \in X$ solving the inclusion $F(x) \ni y$. 

First, one can apply the criterion for single-valued mappings to the restriction of the canonical projection from $X \times Y$ onto $Y$ to the graph of  $F$, that is, the assignment ${\rm gph}\, F\ni(u,v)\longmapsto v\in Y$, e.g., see \cite[Proposition 3]{i2000}; arriving at the statement in the setting $(S_1)$. Equivalently, one can apply Theorem~\ref{Ekeland} to the function $ \gph F \ni (u,v) \longmapsto \varphi (u,v) := \varrho(y, v)$.   
In \cite{cf15}, it was demonstrated  that, in this way, one gets short and transparent proofs of well-known regularity statements by using \emph{three} ingredients only: the definition of the approximating object (e.g., the one-sided directional derivative, the graphical derivative, the contingent variation, etc.), the triangle inequality, and the criterion.      

Second, in the setting $(S_2)$, one could consider the function $X \ni u \longmapsto \varphi(u) := \dist(y, F(u))$ which, unfortunately, is not lower semi-continuous in general. To cure this drawback, the authors work with a lower semi-continuous envelope of $\varphi$ which brings an extra intricacy to the proof.

Third, a much less explored choice is the distance function generated by the graph of $F$, that is, the function $X \ni u \longmapsto \varphi(u) := \dist\big((u,y), \gph F \big)$, which is Lipschitz continuous with the constant $1$. This originates in the notion of graphical regularity introduced by L. Thibault \cite{Thibault99},  which is equivalent to the property of metric regularity. The strength of this approach can be seen from \cite{FIR2018} by M. Fabian, A.D. Ioffe, and J. Revalski, where a separable reduction of local metric regularity is established. We demonstrate that, in this way, one can avoid conditions based on local slopes as well as lower semi-continuous envelopes when proving stability of metric regularity under additive perturbations in the setting $(S_2)$.    
   
\medskip 

To conclude general comments,  we express a personal authors' impression. Whenever one reads a~proof relying on slope-based conditions, there is a~place in the proof, where one can \emph{directly apply a suitable regularity criterion} to get a shorter transparent proof, that is, the following detour does not seem to be necessary: \emph{Ekeland's variational principle $\implies$ the error bound property $\implies$ the slope-based condition for regularity $\implies$ a particular regularity statement}. 

\medskip 

This note is organized as follows. In Section~\ref{scGRC}, we recall various versions of regularity on a fixed set and provide sufficient conditions for these properties in settings $(S_1)$ and $(S_2)$, respectively. In Section~\ref{scPS}, we apply  general criteria to show the stability of regularity under additive perturbations by particular single-valued as well as set-valued mappings. In the last section, we briefly discuss several possible generalizations   and provide an outlook for future work.
     
\medskip 

\paragraph{\bf Notation and terminology.}  
We use the convention that $ \inf \emptyset: = \infty$ and $\sup \emptyset: = 0$, because we work with non-negative quantities; and that $0 \cdot \infty = 1$. If a set is
a singleton we identify it with its only element, that is, we write $a$ instead of $\{ a \}$. By  ${f: X\to Y}$ we denote a single-valued mapping  $f$ acting from a set $X$ into another set $Y$;
while  ${F:X \rightrightarrows Y}$ denotes a~mapping from $X$ into $Y$ which may be set-valued, i.e., the correspondence  determined by a subset of $X \times Y$ called the  \emph{graph} of $F$ and denoted by  $\mbox{\rm gph} \,  F$; the \emph{value} of $F$ at an $x \in X$ is the set $F(x) := \{y \in Y: \, (x,y) \in \mbox{\rm gph} \,  F\}$; the {\it domain} of $F$ is the set $\mbox{\rm dom} \,  F:=\{ x \in X:\;  F(x) \neq \emptyset\}$; and the \emph{inverse} of $F$ is the mapping $Y \ni y \longmapsto \{x \in X : \ y\in F(x)\}=: F^{-1}(y) \subset X$;
thus $F^{-1} :Y\rightrightarrows X$ maps $\rge  F$ into $\dom  F$. 

In a metric space $(X,d)$, by $B_X[x,r]$ and $B_X(x,r)$ we denote the \emph{closed ball} and the \emph{open ball} centered at an $x \in X$ with a radius $r \in (0, \infty]$; in particular,   $B_X[x,\infty] = B_X(x,\infty)=X$; 
the \emph{distance from a point} $x \in X$ to a set $C \subset X$ is $\mbox{\rm dist} \, (x,C):= \inf \{d(x,u): \ u \in C \}$; 
and the \emph{diameter} of $C$ is defined as 
$\diam C = \sup\{d(u,x): u, x \in C \}$. The set $C$ is said to be closed [complete] around an $x \in C$ provided that there is an $r > 0$ such that the set $C \cap \ball_X[x,r]$ is closed [complete].

We  consider only the vector (linear) spaces $X$ the scalar field of which is $\mathbb{R}$. Let $\ell > 0$. 
A set-valued mapping $F$ from a metric space $(X,d)$ into a linear metric space $(Y, \varrho)$ is said to have \emph{Aubin property on a set} $U \times V \subset X \times Y$ with the constant $\ell$ provided that 
 $F(x)\cap V \subset F(u) + \ball_Y[0,\ell \, d(x, u)]$ for each $x$, $u \in U$; to have \emph{Aubin property around a point} $(\bar{x}, \bar{y}) \in \gph F$ with the constant $\ell$ provided that  there are $\alpha > 0$ and $\beta > 0$ such that $F$ has Aubin property on $\ball_X(\bar{x},\alpha) \times  \ball_Y(\bar{y},\beta)$ with the constant $\ell$; and to be \emph{Hausdorff-Lipschitz} on a set $C \subset X$ with the constant $\ell$  provided that 
 $F$ has Aubin property on $C \times Y$  with the constant $\ell$.  

\section{General Regularity Criteria} \label{scGRC}

\subsection{Regularity on a fixed set}  
There are several distinct concepts of regularity on a prescribed set. Given a set $W \subset X \times Y$, we denote by $W_X$ and $W_Y$ the projections of $W$ onto $X$ and $Y$, respectively, defined by
$W_X := \{ x \in X: \ (x,y) \in W \mbox{ for some } y \in Y \}$  and $W_Y := \{y \in Y: \ (x,y) \in W \mbox{ for some } x \in X \}$;
and for a fixed $y \in Y$, we set $  W_{\cdot, y} := \{x \in X: \ (x,y) \in W \}$. Let us start with two basic concepts in the spirit of \cite[p. 328]{DontchevBook}.

\begin{definition} \label{defReg}
  Let $(X,d)$ and $(Y, \varrho)$ be metric spaces. Consider a $\kappa > 0$, non-empty sets $V\subset Y$ and $W \subset  X \times Y$, and a set-valued mapping $F:X\tto Y$. The mapping $F$ is said to be 
  \begin{itemize}
       \item[$(I)$]  \emph{$V-$restrictedly regular on $W$ with the constant $\kappa$} provided that 
       \begin{eqnarray} \label{defMRA}
		\dist \big(x,F^{-1}(y) \big)\leq  \kappa \dist \big(y, F(x) \cap V\big)
		\quad \mbox{for each} \quad (x,y) \in W;
	\end{eqnarray}
	    \item[$(II)$]  \emph{regular on $W$ with the constant $\kappa$} provided that 
       \begin{eqnarray} \label{defMRB}
		\dist \big(x,F^{-1}(y) \big)\leq  \kappa \dist \big(y, F(x)\big)
		\quad \mbox{for each} \quad (x,y) \in W.
	\end{eqnarray}
  \end{itemize}
  Given a point $(\bar{x}, \bar{z}) \in \gph F$, the mapping $F$ is said to be \emph{regular around $(\bar{x}, \bar{z})$ with the constant $\kappa$}  provided that there is a neighborhood $W$ of  $(\bar{x}, \bar{z})$ in $X \times Y$ such that $F$ is regular on $W$ with the constant $\kappa$.  
\end{definition}

\begin{remark} \rm Let us emphasize that one has to impose additional assumptions on $V$, $W$, and  $F$  in order to obtain non-trivial properties. For example, (I) holds for \emph{each} $F$ with $F(x) \cap V = \emptyset$ for each $x \in W_X$; and the same is true in  (II) when  $\dom F \cap W_X = \emptyset$. Obviously, if $F$ is regular on $W$ with the constant $\kappa$, then $F$ is $V-$restrictedly regular on $W$ with the constant $\kappa$ for each non-empty $V \subset Y$. On the other hand, $F$ is $V-$restrictedly regular on $W$ with the constant $\kappa$ provided that the mapping $\widetilde{F}(x) := F(x)\cap V$, $x \in X$, is regular on $W$ with the constant $\kappa$.  
\end{remark}

There is a substantial difference between properties (I) and (II) as the following easy but leading example shows.
\begin{example} \label{exFregIandII}  \rm 
 Let $(X,d):= (Y,\varrho) := (\mathbb{R},\vert \cdot\vert)$ and $(\bar{x},\bar{z}):=(0,0)$. Consider a mapping $\mathbb{R}\ni x \longmapsto F(x) := \lbrace x,-1 \rbrace\subset \mathbb{R}$. Given $\alpha > 0$ and $\beta>0$, we investigate regularity of $F$ on a neighborhood $ W(\alpha, \beta) := (-\alpha, \alpha) \times (-\beta, \beta)$ of $(\bar{x}, \bar{z}) \in \gph F$.

First, pick a $\beta > 1$. Then there are no $\alpha > 0$ and $\kappa > 0$ such that $F$ is regular on $W(\alpha, \beta)$ with the constant $\kappa$. Assume, on the contrary, that there are such constants $\alpha$ and $\kappa$. Pick a $t > 0$ such that $y := -1 -t \in (-\beta, \beta)$ and $x := \kappa t \in  (-\alpha, \alpha)$. Then 
 $$
   \dist(x, F^{-1}(y)) = |x-y| = \kappa  t + t +1  >  \kappa t =   \kappa \,  |y + 1| = \kappa \dist(y,F(x)),  
$$ 
a contradiction.  Similarly,  given an $r>1$, there are no $\alpha > 0$ and $\kappa >0 $ such that $F$ is $(-r,r)$-restrictedly regular on $W(\alpha, \beta)$ with the constant $\kappa$. On the other hand, for any $\alpha > 0$ and any $r\in(0,1)$,  the mapping $F$ is $(-r,r)$-restrictedly regular on $W(\alpha, \beta)$ with the constant $1$.
 
Second, pick a $\beta \in (0, 1/2)$ and let $\alpha := 1/2 - \beta > 0$.  Given an arbitrary $(x,y) \in W(\alpha, \beta)$, we have 
$ \dist(x, F^{-1}(y)) = |x-y| =  \dist(y,F(x))$. Hence $F$ is regular on $W(\alpha, \beta)$ with the constant $1$; in particular, $F$ is regular around $(\bar{x},\bar{z})$ with the constant $1$. 

Finally, fix a  $\beta\in (1/2,1)$. Let $\alpha>0$ be arbitrary. Then 
\begin{eqnarray}
	\label{eqRerModulusEstimated}
	\kappa=\kappa(\alpha, \beta):=\sup_{(x,y)\in W(\alpha, \beta)}\tfrac{\vert x-y\vert}{  y+1}\geq \sup_{y\in (-\beta, \beta)}\tfrac{\vert y\vert}{  y+1}=\tfrac{\beta}{1-\beta} > 1.
\end{eqnarray}
 Given a point  $(x,y)\in W(\alpha, \beta)$, we have $\dist(y,F(x))=  \min\lbrace y+1,\vert x-y\vert\rbrace$; thus 
\begin{eqnarray*}
	\dist(x,F^{-1}(y))=\vert x-y\vert\leq \kappa   \min\lbrace y+1,\vert x-y\vert\rbrace=\kappa\dist(y,F(x)).
\end{eqnarray*}
 Hence $F$ is regular on $W(\alpha, \beta)$. By \eqref{eqRerModulusEstimated}, we get the constant of regularity tends to $\infty$ as $\beta \uparrow 1$ without regard to the value of $\alpha$, that is, $F$ is losing regularity. On the other hand, for any $\alpha > 0$ and any $r\in(0,1)$,  the mapping $F$ is $(-r,r)$-restrictedly regular on $W(\alpha, \beta)$ with the constant $1$.
\end{example}

In \cite{Ioffe10}, A.D. Ioffe defined an abstract version of regularity involving not only the set $W$ as above but also containing a reference point and a radius of a ball around this point as extra parameters; unifying in this way almost all regularity concepts. However, we prefer a top-to-bottom approach, that is, to introduce these extra parameters at the time when (semi-)local statements are considered, in order to keep the presentation as transparent as possible.  The property (II) covers, in particular,  the so-called $\gamma$-regularity and Milyutin regularity.
 
\begin{definition} \label{defRegBtypes}
   Let $(X,d)$ and $(Y, \varrho)$ be metric spaces. Consider constants $\delta>0$ and $\kappa>0$, a non-empty set $W \subset  X \times Y$, a set-valued mapping $F:X\tto Y$, and a function $\gamma: X \to [0, \infty]$ being positive on $\dom F \cap W_X \neq \emptyset$. The mapping $F$ is said to be 
  \begin{itemize} 
     \item[$(a)$] \emph{$\gamma$-regular on $W$ with the constant $\kappa$} provided that 
      $F$ is regular with the constant $\kappa$ on the set 
	$\big\{(x,y) \in W: \ \kappa \dist \big( y,F(x) \big) < \gamma(x) \big\}$;
  	\item[$(b)$] \emph{$\gamma$-regular on $W$ with the constants $\delta$ and $\kappa$} provided that 
      $F$ is regular with the constant $\kappa$ on the set 
     $\big\{(x,y) \in W: \  \dist \big( y,F(x) \big) < \delta \gamma(x) \big\}$.
	\end{itemize}
	If $\gamma(u) := \dist(u, X \setminus W_X)$, $u \in X$, either in $(a)$ or (b); then we speak about \emph{Milyutin regularity}.
\end{definition}

\begin{remark} \rm 
Consider non-empty sets $U \subset X$ and $V \subset Y$, a constant $\kappa >0$,   and a mapping $F: X \rightrightarrows Y$. 
Put $W := U \times V$, that is, $W_X = U$ and $W_Y = V$.  Inequality \eqref{defMRA} appears in the definition of metric regularity of $F$ on $U$ for $V$ \cite[p. 328]{DontchevBook}, while  \eqref{defMRB} corresponds to the alternative definition therein; however these properties request, in addition, that the set $\gph F \cap W$ is closed in $X \times Y$. Further, (a) reduces to $\gamma$-metric regularity on $(U,V)$ \cite[Definition 2.22]{IoffeBook}; 
while Milyutin case corresponds to \cite[Definition 2.28]{IoffeBook}. If $U := \ball_X(\bar{x}, r)$, then $\gamma(u) := \dist(u, X \setminus U) = r - d(u, \bar{x})$ for each $u \in U$. Any $F$ having property in (a) has property in (b) for any  $\delta \in (0, 1/\kappa]$; and the same is true for Milyutin regularities.  
\end{remark}

Note that, setting $\gamma(x) := \infty$, $x \in W_X$, property in (a) is nothing else but regularity of $F$ on $W$ with the constant $\kappa$ from (II). However, we prefer to consider (a) as a particular case of (II) because the reversed process would give weaker sufficient conditions in the following subsections. 

\subsection{Mappings with a complete graph} \label{sscGRCcg} A stronger  concept of regularity  than (II) in the setting $(S_1)$ was studied in \cite{Ioffe10}, see also Section~\ref{scCR}. Therefore we focus on the restricted regularity and  start with single-valued mappings.

\begin{proposition} \label{propSingle}
	Let $(X,d)$ be a complete metric space and $(Y, \varrho)$ be a metric space. 
	Consider a $\kappa > 0$, non-empty sets $V \subset Y$ and $W \subset X \times Y$,  and a~continuous single-valued mapping $g:X \to Y$, defined on whole $X$ and such that $g(W_X) \cap V \neq \emptyset$.
	Then $g$ is $V$-restrictedly regular on $W$ with the constant $\kappa$ provided that for each
	$u \in X$ and each $y\in W_Y$ such that $0<  \kappa  \varrho(y,g(u)) \leq \kappa  \varrho(y,g(x)) - d(u,x)$ for some $x \in W_{\cdot,y}$ with $g(x)\in V$, there is a $\hat{u} \in X$, better than $u$, satisfying
	\begin{eqnarray}
		\label{estimateEQMS3M}
		 \kappa\varrho(y,g(\hat{u})) < \kappa \varrho(y,g(u)) - d(\hat{u},u).
	\end{eqnarray}
\end{proposition}
\begin{proof}
	Pick any $(x,y) \in W$. If $g(x) = y$ or $g(x) \notin V$, then we are done. Assume that $y \neq g(x) \in V$. 
	Theorem \ref{Ekeland}, with (a continuous) $\varphi:= \kappa \varrho(y, g(\cdot))$, yields a~$u\in X$ such that
	\begin{equation}
		\label{EstimateEXMS1M}
		\kappa \varrho(y,g(u)) + d(u,x)\leq \kappa\varrho(y,g(x))
	\end{equation}
	and 
	\begin{equation}
		\label{estimateEQMS5M}
		\kappa\varrho(y,g(u)) \leq \kappa\varrho(y,g(\hat{u}))+ d(\hat{u},u)\quad\text{whenever}\quad \hat{u}\in X.
	\end{equation}
	By \eqref{EstimateEXMS1M}, we have $d(x,u)\leq \kappa \varrho (y,g(x))$.  To get the inequality in \eqref{defMRA} with $F := g$, it suffices to show that $u \in g^{-1}(y)$. Assume, on the contrary, that $g(u) \neq y$.	
 Remembering that $x \in W_{\cdot, y}$ is such that  $g(x)\in V$ and that $y \in W_Y$, using  \eqref{EstimateEXMS1M}, we infer that 
	$$
	  0 < \kappa\varrho(y,g(u))  \leq \kappa\varrho(y,g(x)) - d(u,x).
	$$
	Hence the assumption containing \eqref{estimateEQMS3M}, yields a $\hat{u}\in X$ such that
	$d(\hat{u}, u)< \kappa\varrho(y,g(u))-\kappa\varrho(y,g(\hat{u}))$. This in combination with \eqref{estimateEQMS5M} implies that
	$$
	 d(\hat{u},u) <\kappa\varrho(y,g(u))-\kappa\varrho(y,g(\hat{u})) \leq  \kappa\varrho(y,g(\hat{u}))+ d(\hat{u},u) - \kappa\varrho(y,g(\hat{u})) = d(u,\hat{u}),
	$$
	a contradiction. Hence $y=g(u)$.
 \end{proof}

Now, we use the ``projection trick'' described in the introduction. 
Let $(X,d)$ and $(Y, \varrho)$ be metric spaces and let $\lambda > 0$. Denote by $\xi$  a~(compatible) metric on $X \times Y$ defined by 
\begin{equation}  \label{defXi}
  \xi((u, v), (x, y)) := \max \{ d(u,x), \lambda \varrho(v,y) \}
  \quad \mbox{for each} \quad  (u,v),  (x, y) \in  X \times Y. 
\end{equation}

\begin{proposition}  \label{propIoffeRegA}	
   Let $(X,d)$  and $(Y, \varrho)$ be a metric spaces.  Consider constants $\lambda > 0$ and  $\kappa>\lambda$, non-empty sets $V \subset Y$ and  $W \subset X \times Y$, and a~set-valued mapping $F:X\rightrightarrows Y$ such that  $\gph F $ is complete and $F(W_X) \cap V \neq \emptyset$. Let us equip $X \times Y$ by the metric $\xi$ in \eqref{defXi}. 		  
    Then $F$ is $V-$restrictedly regular on $W$ with the constant $\kappa$ provided that   for each $(u,v) \in \gph F$ and each $y\in W_Y$ such that $0<\kappa \varrho(y,v) \leq \kappa \varrho(y,\tilde{v})- \xi((u,v),(x,\tilde{v}))$ for some $x \in W_{\cdot, y}$ and $\tilde{v} \in F(x) \cap V$, there is a pair $(\hat{u},\hat{v})\in\gph	F $, better than $(u,v)$, satisfying
	\begin{equation}
		\label{LinOpenSetEstimM}
		\kappa\varrho(y,\hat{v}) <\kappa\varrho(y,v) - \xi((\hat{u},\hat{v}), ( u,v)). 
	\end{equation}
\end{proposition}

\begin{proof} We follow \cite[Proposition 7]{Ioffe10} and the proof of it. Let $\widetilde{X}:=  \gph F$ and 
 $\Omega := \{ ((x,v),y) \in \widetilde{X} \times Y:  \ v \in V \mbox{ and } (x,y) \in W \}$.   Then $(\widetilde{X}, \xi)$ is a~complete metric space, $\Omega_{\widetilde{X}} \subset \widetilde{X} \cap (W_X \times V)$, $\Omega_Y \subset W_Y$, and $\Omega_{\cdot, y} \subset W_{\cdot, y}\times V$ for each $y \in Y$.	
 
 	We apply Proposition~\ref{propSingle} with $(X, d):=(\widetilde{X}, \xi)$ and a continuous mapping $\widetilde{X}\ni (x,v) \longmapsto g(x,v) := v$.   	
	Note that $V \supset g(\Omega_{\widetilde{X}}) \neq \emptyset$. Fix any $(u,v)\in \widetilde{X}$ and any $y\in \Omega_Y$ such that  $0< \kappa \varrho (y,g(u,v)) \leq\kappa \varrho (y,g(x,\tilde{v})) -\xi ((u,v), (x,\tilde{v}))$ for a pair $(x,\tilde{v}) \in \Omega_{\cdot,y}$. Then  $(u,v) \in \gph F$, $y \in W_Y$, $x \in W_{\cdot,y}$,  $g(u,v)=v$, and $g(x,\tilde{v})=\tilde{v} \in F(x) \cap V$.  Find a pair $(\hat{u},\hat{v})\in\gph F$ satisfying \eqref{LinOpenSetEstimM}, which is \eqref{estimateEQMS3M} with $d:=\xi$ for our $g$. 
	Applying  Proposition~\ref{propSingle}, we get that 
	$$
	\dist \big((x,v),g^{-1}(y)\big) \leq \kappa \varrho(y, g(x,v))
       \quad \mbox{whenever} \quad ((x,v),y) \in \Omega. 
	$$
	Fix arbitrary $(x,y) \in W$. Pick a $v \in F(x) \cap  V$ (if there is any). Then $((x,v),y) \in \Omega$ and   $\dist \big(x,F^{-1}(y)\big)\leq \dist \big((x,v),g^{-1}(y)\big)$. The last displayed formula yields that
	$
		\dist \big(x,F^{-1}(y)\big)\leq \dist \big((x,v),g^{-1}(y)\big)\leq \kappa \varrho(y,g(x,v)) = \kappa\varrho(y,v).
	$
	As $v \in F(x) \cap V$ is arbitrary, the proof is finished. 
\end{proof}

In the literature, the second inequality in the constraint on $(u,v)$ and $y$ is commonly omitted, that is, only $v \neq y$ is used. However, this immediately excludes the possibility to get effective statements for $W_X = \{\bar{x}\}$, that is, for the property of \emph{semiregularity}, see \cite{CFK2018} for a detailed treatment and a list of relevant references. 

\begin{corollary}  \label{corIoffeSemiRegA}	
   Let $(X,d)$  and $(Y, \varrho)$ be a metric spaces.  Consider constants $\lambda > 0$ and  $\kappa > \lambda$, non-empty sets $\Gamma$, $\Lambda \subset Y$,  and a~set-valued mapping $F:X\rightrightarrows Y$ such that  $\gph F $ is complete and $ F(\bar{x}) \supset \Lambda$. Let us equip $X \times Y$ by the metric $\xi$ in \eqref{defXi}. 		  
    Assume that  for each $(u,v) \in \gph F$ and each $y\in \Gamma$ such that $0<\kappa \varrho(y,v) \leq  \kappa \varrho(y,\tilde{v})- \xi((u,v),(\bar{x},\tilde{v}))$ for some $\tilde{v} \in \Lambda$, there is a pair $(\hat{u},\hat{v})\in\gph	F $, better than $(u,v)$, satisfying \eqref{LinOpenSetEstimM}. Then   
    \begin{equation} \label{defSemireg}
    	\dist \big(\bar{x},F^{-1}(y)\big)\leq  \kappa \dist ( y, \Lambda )
    	\quad \mbox{for each} \quad y\in \Gamma.
    \end{equation}		
\end{corollary}
\begin{proof} 
  Let $W := \{\bar{x}\} \times \Gamma$ and $V := \Lambda$.  Then $W_Y = \Gamma$ and  $W_{\cdot, y} = \{\bar{x}\}$ for each $y \in W_Y$. Proposition~\ref{propIoffeRegA} yields the conclusion.  
\end{proof}

The set $\Gamma$ represents the set of the admissible right-hand sides $y$ for which we want to solve the inclusion $F(x) \ni y$; while $\Lambda$, ranging between $\{\bar{y}\}$ and $F(\bar{x})$, determines the strength of the estimates on the distance from $\bar{x}$ to the solution set $F^{-1}(y)$. Taking $\Gamma := \ball_Y(\bar{y}, r/\kappa)$ and $\Lambda := \{\bar{y}\}$, with fixed $r > 0$ and $\bar{y} \in F(\bar{x})$, we arrive at \cite[Proposition 4.2(i)]{CFK2018}.

\subsection{Mappings with a complete domain and a closed graph} \label{sscGRCcdcg}

Let $(X,d)$ and $(Y, \varrho)$ be metric spaces and let $\kappa > 0$. Denote by $\omega$  a (compatible) metric on $X \times Y$ defined by 
\begin{equation}  \label{defOmega}
  \omega((u, v), (x, y)) := d(u,x) + \kappa \varrho(v,y)
  \quad \mbox{for each} \quad  (u,v),  (x, y) \in  X \times Y. 
\end{equation}

We start with Ioffe-type conditions for regularity in (II) in the setting $(S_2)$. 
   
\begin{proposition}  \label{propIoffeGraph}	Let $(X,d)$ be a complete metric space and $(Y, \varrho)$ be a metric space. 
		Consider constants  $\kappa > 0$ and $\hat{\kappa} \geq \kappa$,  a set $W \subset X \times Y$,
		and  a~set-valued mapping $F:X\rightrightarrows Y$, with a closed graph, such that $\dom F \cap W_X \neq \emptyset$.  Let us equip $X \times Y$ by the metric $\omega$ in \eqref{defOmega}. 		
	Then $F$ is regular on $W$ with the constant $\hat{\kappa}$  provided that  for  each $u \in X$ and each $y \in W_Y$ such that  
	\begin{equation}
		\label{eqBad}
	  0<  \dist \big( (u,y), \gph F \big) \leq \kappa \dist\big(y, F(x)\big) - \frac{\kappa}{\hat{\kappa}}  \, d(u,x)
	\end{equation}
  for some $x \in W_{\cdot, y}$,   there is a  $\hat{u} \in X$, better than $u$, satisfying
	\begin{equation}
		\label{eqBetter}
		\dist \big((\hat{u},y), \gph F \big) < \dist \big((u,y), \gph F \big) - \frac{\kappa}{\hat{\kappa}}  \,  d(\hat{u},u). 
	\end{equation}
\end{proposition}
\begin{proof}
  Put $\alpha := \kappa/\hat{\kappa} \in (0,1]$. Fix any $(x, y) \in W$. If $ y \in F(x)$,  then we are done. Assume that $y \notin F(x)$. The function 
 $
   (X,d) \ni u  \longmapsto \varphi(u):=\dist \big( (u,y), \gph F \big)
 $ 
 is Lipschitz continuous with the constant $1$  and   $\varphi(u) \leq \kappa \dist\big( y,F(u) \big)$ for each $u  \in X$. Also $\varphi(x) > 0$ since the graph of $F$ is closed.   Theorem \ref{Ekeland}, with $d:= \alpha \, d$,  yields a $u\in X$ such that
	\begin{equation}
		\label{eqEVP01}
		\varphi(u) + \alpha \,  d(u,x)\leq \varphi(x)
	\end{equation}
	and
	\begin{equation}
		\label{eqEVP02}
		 \varphi(u) \leq \varphi(\hat{u}) + \alpha\,  d(\hat{u},u)\quad\text{whenever}\quad \hat{u}\in X.
	\end{equation}
	By \eqref{eqEVP01}, we have
	\begin{equation}
		\label{eqEVP01d}
		\dist \big( (u,y), \gph F \big) +	\alpha  \,  d (u,x)\leq \varphi(x) \leq  \kappa \dist \big( y,F(x) \big). 
	\end{equation}
	In particular, 
	\begin{equation}
		\label{eqR00}
	  d(x,u) \leq  \hbox{$\frac{\kappa}{\alpha}$}  \dist\big( y,F(x) \big)  = \hat{\kappa} \dist\big( y,F(x) \big). 
	\end{equation}
	To get the inequality in \eqref{defMRB} with $\kappa := \hat{\kappa}$, it suffices to show that $u  \in F^{-1}(y)$. Assume, on the contrary, that $y \notin F(u)$.   Since $(x,y) \in W$, we have $y \in W_Y$ and $x \in W_{\cdot, y}$.
	As the graph of $F$ is closed, \eqref{eqEVP01d}  yields that    
	\begin{eqnarray*}
	     0 & < & \dist \big( (u,y), \gph F \big) \leq   \kappa \dist\big( y,F(x) \big) - \alpha \,  d(u,x), 
	\end{eqnarray*}
  	which is \eqref{eqBad}. Thus  there is  a $\hat{u}\in X$ such that
	$ \alpha \, d(\hat{u},u)< \varphi(u) - \varphi(\hat{u})$.
	This in combination with  \eqref{eqEVP02} implies that
	$$
    	\alpha \, d(\hat{u},u)  <  \varphi(u) - \varphi(\hat{u})   \leq  \varphi (\hat{u}) + \alpha \, d(\hat{u},u) - \varphi(\hat{u})  = \alpha \, d(\hat{u},u),
	$$
	a contradiction. Hence $y\in F(u)$. 
\end{proof}
 
The above statement is a quantitative non-local version of the sufficiency part of \cite[Theorem~2]{FIR2018}, where the second inequality in \eqref{eqBad} is omitted, and covers $\gamma$-regularity and, in particular,  Milyutin regularity; while the corresponding statement in the setting $(S_1)$, with $W := U \times V$, is \cite[Theorem 3.1]{Ioffe11}.   
 
\begin{corollary}  \label{corIoffeGraph1}
  Let $X$, $Y$, $\kappa$, $\hat{\kappa}$, $W$, $F$, and $\omega$ be as in Proposition~\ref{propIoffeGraph}. Consider a function $\gamma : X \to [0, \infty]$ which is positive  on $\dom F \cap W_X$.     Then $F$ is $\gamma$-regular on $W$ with the constant $\hat{\kappa}$  provided that  for  each $u \in X$ and each $y \in W_Y$ such that  $0< \hat{\kappa} \dist \big( (u,y), \gph F \big) < \kappa \big(\gamma(x) - d(u,x)\big)$  for some $x \in W_{\cdot, y}$,   there is a  $\hat{u} \in X$, better than $u$, satisfying \eqref{eqBetter}.
\end{corollary}
\begin{proof}  Put $\alpha := \kappa/\hat{\kappa} \in (0,1]$. We apply  Proposition~\ref{propIoffeGraph} with $W$ replaced by 
	$$
     \Omega := \big\{(x,y) \in W: \ \hat{\kappa} \dist \big( y,F(x) \big) < \gamma(x) \big\}.
	$$
    Pick any $u \in X$ and any $y \in   \Omega_Y \subset W_Y$ such that \eqref{eqBad} holds true for some $x \in \Omega_{\cdot, y} \subset W_{\cdot, y}$. 
   Then  \begin{eqnarray*}
	     0 & < & \dist \big( (u,y), \gph F \big) \leq   \kappa \dist\big( y,F(x) \big) - \alpha \,  d(u,x) \\
	     & = & \alpha  \, \big( \hat{\kappa} \dist \big( y,F(x) \big) - \, d(u,x) \big)  <    \alpha \, \big( \gamma(x) -   d(u,x) \big).
	\end{eqnarray*}
   Thus there is a $\hat{u}\in X$  satisfying \eqref{eqBetter}. Proposition~\ref{propIoffeGraph} says that $F$ is regular on $\Omega$ with the constant $\hat{\kappa}$, that is,  $F$ is $\gamma$-regular on $W$ with the constant $\hat{\kappa}$.
	\end{proof}

In case of Milyutin regularity,  the above statement provides a conclusion for the \emph{interior} points of $W_X$ only. Consider a set $U \subset X$ with a non-empty interior, a non-empty set $V \subset Y$, and a point $(\bar{x},\bar{y}) \in X \times Y$. The choice  $W := U \times \{\bar{y}\}$ is allowed and provides conditions for the property called \emph{subregularity} while the choice $W := \{\bar{x}\} \times V$ yields  a void statement, that is, the property of \emph{semiregularity} is not covered. On the other hand, Proposition~\ref{propIoffeGraph} gives an analogue of Corollary~\ref{corIoffeSemiRegA}. 

\begin{corollary}  \label{corIoffeGraphSemireg}	Let $(X,d)$ be a complete metric space and $(Y, \varrho)$ be a metric space. 
		Consider constants  $\kappa > 0$ and $\hat{\kappa} \geq \kappa$,  a non-empty set $\Gamma \subset Y$, a point $\bar{x} \in X$,  and a~set-valued mapping $F:X\rightrightarrows Y$, with a closed graph, such that $\dom F \ni \bar{x}$.  Let us equip $X \times Y$ by the metric $\omega$ in \eqref{defOmega}. 		
Assume 	that for  each $u \in X$ and each $y \in \Gamma$ such that   
	$	  0 <  \dist \big( (u,y), \gph F \big) \leq \kappa \dist\big(y, F(\bar{x})\big) - \frac{\kappa}{\hat{\kappa}} \, d(u, \bar{x})$,
	there is a  $\hat{u} \in X$, better than $u$, satisfying \eqref{eqBetter}. 
	Then \eqref{defSemireg}, with $\kappa := \hat{\kappa}$ and $\Lambda := F(\bar{x})$, holds true.
\end{corollary}
\begin{proof} 
  Let $W := \{\bar{x}\} \times \Gamma$.  Then $W_Y = \Gamma$ and  $W_{\cdot, y} = \{\bar{x}\}$ for each $y \in W_Y$. Proposition~\ref{propIoffeGraph} yields the conclusion.  
\end{proof}

\section{Perturbation Stability} \label{scPS}

\subsection{Mappings with a complete graph} \label{sscPScg}
In \cite{Ioffe10}, there is studied the stability of a  stronger version of regularity than  (II), cf. Section~\ref{scCR}, under single-valued additive perturbations in the setting $(S_1)$. Therefore we focus on the stability of the regularity in (I) under particular set-valued additive perturbations.  We start with a semi-local result having several seemingly different reincarnations.   In the proof, we shall employ a truncated version of the so-called \emph{epigraphical multifunction} from \cite{HNT2014} which allows to avoid the assumption that the sum under consideration has a closed graph.   

	\begin{theorem} \label{thmSUSVP_S1}
	 Let $(X,d)$ be a metric space and $(Y,\varrho)$ be a linear metric space with a translation-invariant metric. Consider positive constants $a$, $b$, $r$, $\kappa$, $\hat{\kappa}$, and $\ell$ such that  $\kappa \ell < 1$ and $\hat{\kappa}>\kappa/(1-\kappa\ell)$,  a $c \in [0,\infty]$,  a point $(\bar{x}, \bar{z}, \bar{w}) \in X \times Y \times Y$, and set-valued mappings $F$, $G:X\tto Y$, with $\bar{z}\in F(\bar{x})$ and $\bar{w} \in G(\bar{x})$. Assume that 
	 $\alpha := a + \hat{\kappa} (b + r)$, $\beta := c + \alpha \ell$, and $\delta := \beta + b + r$
	are such that $\gph F\cap\big(\ball_X[\bar{x},\alpha]\times \ball_Y[\bar{z},  \delta ]\big)$ and $\gph G\cap\big(\ball_X[\bar{x},\alpha]\times \ball_Y[\bar{w},\beta]\big)$ are complete and that the following conditions hold true: 
	\begin{itemize}
		\item[$(i)$] $F$ is $\ball_Y[\bar{z},\delta]-$restrictedly regular on $ \ball_{{X}}(\bar{x}, \alpha) \times \ball_{{Y}}(\bar{z},\delta) $  with the constant $\kappa$;
		\item[$(ii)$]  $G$ has Aubin property on $\ball_{{X}}(\bar{x}, \alpha) \times \ball_Y(\bar{w}, \beta)$ with the constant $\ell$; 
		\item[$(iii)$] for each $x \in \ball_{{X}}(\bar{x}, a)$ and each $v \in (F+G)(x) \cap \ball_{{Y}}[\bar{z} + \bar{w}, b + r]$ there is a $w \in G(x) \cap \ball_Y(\bar{w}, a\ell + c )$ such that $v \in F(x) + w$.
	\end{itemize}
    Then  for each $\bar{q} \in \ball_Y\big[\bar{w}, \min\{r,b\}\big]$ the sum $F+G$ is  $\ball_Y[\bar{z} + \bar{q},r]-$restrictedly regular on $  \ball_X(\bar{x},a) \times \ball_Y(\bar{z} + \bar{q},b)$ with the constant $\hat{\kappa}$. 
\end{theorem}
	\begin{proof} Fix any $\bar{q}$ as above.
		 Find a $\tilde{\kappa} \in(\kappa, 1/\ell)$  such that 
		 $\hat{\kappa} > \tilde{\kappa}/(1-\tilde{\kappa} \ell)> \kappa/(1-\kappa \ell)$. 	Let $\lambda:=\tilde{\kappa}/(1+\tilde{\kappa}\ell)$. Define a metric $\tilde{d}$ on  $\widetilde{X}:={X} \times Y$  by
		\begin{eqnarray*}
			\tilde{d}((u,w),(\hat{u},\hat{w})):=\max\lbrace d(u,\hat{u}),  \varrho(w,\hat{w})/\ell\rbrace\quad\text{for each}\quad (u,w),(\hat{u},\hat{w})\in \widetilde{X};
		\end{eqnarray*}
		sets  $W:=\ball_{X}(\bar{x},a) \times \ball_Y(\bar{w},a \ell + c) \times \ball_Y(\bar{z}+\bar{q},b)$ and $V :=\ball_Y[\bar{z} + \bar{q},r]$; and a~mapping $\mathcal{E}: \widetilde{X} \tto Y$ by 		
		$$
		\mathcal{E}(x,w):= F(x)\cap \ball_{{Y}}[\bar{z}, \delta]    + w     \mbox{ if }  x \in \ball_X[\bar{x}, \alpha] \mbox{ and } w \in G(x)\cap\ball_Y[\bar{w}, \beta];
		$$
		and $\mathcal{E}(x,w)=\emptyset$ otherwise. Then $\mathcal{E}$ has a complete graph.   We check the assumptions of Proposition \ref{propIoffeRegA}, with $F:=\mathcal{E}$, $(X,d):=(\widetilde{X},\tilde{d})$,  and $\kappa:=¨\tilde{\kappa}/(1-\tilde{\kappa} \ell)$, where $\xi$ is defined by \eqref{defXi} with $d := \tilde{d}$ and $X := \widetilde{X}$.  Note that $\lambda < \tilde{\kappa} < \tilde{\kappa}/(1 - \tilde{\kappa}\ell)$ and $W_{\widetilde{X}} = \ball_{X}(\bar{x}, a) \times \ball_Y(\bar{w}, a \ell + c)$. As $\varrho(\bar{q}, \bar{w}) \leq r$,  we infer that $\bar{z} + \bar{w} \in \mathcal{E} (W_{\widetilde{X}}) \cap V \neq \emptyset$.
		Fix arbitrary $(u, w, v)\in \gph \mathcal{E}$  and $y\in\ball_{{Y}}(\bar{z}+\bar{q},b)$ such that 
        \begin{eqnarray}
        	\label{eqFixedEstimated}
		\quad 0<\tilde{\kappa}\varrho(y,v) \leq \tilde{\kappa}\varrho(y,\tilde{v})-(1- \tilde{\kappa} \ell ) \max\lbrace d(u,x), \varrho(w,\tilde{w})/\ell,\lambda\varrho(v,\tilde{v}) \rbrace
	    \end{eqnarray}
 for some $(x,\tilde{w})\in \ball_{X}(\bar{x}, a) \times \ball_Y(\bar{w}, a \ell + c)$ and $\tilde{v}\in\mathcal{E}(x,\tilde{w})\cap \ball_{{Y}}[\bar{z}+\bar{q},r]$.
		We have to find a triplet $(\hat{u},\hat{w}, \hat{v})\in \gph \mathcal{E}$ such that 
		\begin{eqnarray}
			\label{eqIoffeEpigraph}
			\tilde{\kappa}\varrho(y,\hat{v})<\tilde{\kappa}\varrho(y,v)- (1-\tilde{\kappa} \ell) \max\lbrace d(\hat{u},u), \varrho(\hat{w},w)/\ell, \lambda \varrho(\hat{v},v) \rbrace.
		\end{eqnarray}
        As $\varrho(y,\tilde{v}) \leq  \varrho(y,\bar{z}+\bar{q}) + \varrho(\bar{z}+\bar{q},\tilde{v}) < b + r$ and $ (1 - \tilde{\kappa} \ell) \hat{\kappa}/\tilde{\kappa} > 1$,    multiplying \eqref{eqFixedEstimated} by $\hat{\kappa}/\tilde{\kappa}$,  we get that   
        \begin{equation} \label{eqDux}
             \max\lbrace d(u,x), \varrho(w,\tilde{w})/\ell \rbrace + \hat{\kappa} \varrho(y,v) \leq  \hat{\kappa} \varrho(y,\tilde{v})<  \hat{\kappa} (b + r).
        \end{equation}        		
        Thus $d(u, \bar{x}) \leq d(u,x) + d(x, \bar{x})  <  \alpha$; 
        $\varrho(w,\bar{w}) \leq \varrho(w,\tilde{w}) + \varrho(\tilde{w},\bar{w}) < \alpha \ell + c$; and
       	\begin{eqnarray*}
			\varrho(y-w,\bar{z})&=& 		\varrho(y-w+\bar{q},\bar{z}+\bar{q})\leq 	\varrho(y-w+\bar{q},y)+	\varrho(y,\bar{z}+\bar{q}) \\
			&  \leq & \varrho(w, \bar{w})+	\varrho(\bar{w}, \bar{q}) + \varrho(y,\bar{z}+\bar{q}) <  \alpha \ell +c + r + b.
		\end{eqnarray*}
		Summarizing, we have that  $u\in\ball_X(\bar{x},\alpha)$, that $w \in G(u) \cap \ball_Y(\bar{w},  \beta)$,  and that $y-w\in \ball_{{Y}}(\bar{z}, \delta)$. As   ${v}-w\in F(u)\cap\ball_Y[\bar{z}, \delta]$ and  $0<\kappa<\tilde{\kappa}$, by (i), there is a~$\hat{u}\in F^{-1}(y- w)$ such that
	\begin{eqnarray}
		\label{eqMetricRegularityEstimated2}
		 \quad d({u},\hat{u})<\tilde{\kappa} \varrho(y-{w},{v}-{w})=\tilde{\kappa} \varrho(y,{v}) \quad \quad \big(<\hat{\kappa} \varrho(y,{v}) \big).
	\end{eqnarray}
		This and \eqref{eqDux}  imply that   $d(\hat{u},\bar{x})\leq d(\hat{u},u)+d(u,x) + d(x,\bar{x}) <  \hat{\kappa} ( b + r) + a$; thus  $\hat{u}\in\ball_X(\bar{x},\alpha)$. 
		By {\rm (ii)} and \eqref{eqMetricRegularityEstimated2}, there is a $\hat{w}\in G(\hat{u})$ such that 
		\begin{equation} \label{eqhatw}
		  \varrho({w},\hat{w})\leq \ell \, d({u},\hat{u}) \quad \quad \big(<\tilde{\kappa} \ell \, \varrho(y,{v}) \big).
		\end{equation}  
		Let $\hat{v}:=y-{w}+\hat{w} \in F(\hat{u})\cap\ball_{{Y}}(\bar{z}, \delta) +\hat{w}$. From \eqref{eqhatw} and \eqref{eqDux}, we get that $ \varrho(\hat{w}, \bar{w}) \leq  \varrho(\hat{w}, w) + \varrho(w, \bar{w}) < \hat{\kappa} \ell \varrho(y,v) + \varrho(w, \tilde{w}) + \varrho(\tilde{w}, \bar{w}) < \hat{\kappa} \ell (b + r) + a \ell + c.$
	    Therefore  $\hat{v} \in \mathcal{E}(\hat{u},\hat{w})$. It remains to prove \eqref{eqIoffeEpigraph}. By \eqref{eqhatw}, we get that  
		\begin{eqnarray*}
			\varrho(v,\hat{v}) &\leq&  \varrho(v,y)+\varrho({w},\hat{w}) < \varrho(v,y)+\tilde{\kappa} \ell \varrho(y,v)=(1+ \tilde{\kappa} \ell )\varrho(v,y). 
		\end{eqnarray*}
		Multiplying this by $\lambda$ and taking into account also \eqref{eqhatw}, we infer that $\tilde{\kappa} \varrho(y,v) > \max\lbrace d(\hat{u},u), \varrho(\hat{w},w)/\ell, \lambda \varrho(\hat{v},v) \rbrace$ and that 
		\begin{eqnarray*}
			\tilde{\kappa} \varrho(y,\hat{v})&=& \tilde{\kappa} \varrho(w,\hat{w}) < \tilde{\kappa} \tilde{\kappa}\ell \varrho(y,v) =  \tilde{\kappa} \varrho(y,v) -(1 - \tilde{\kappa}\ell) \tilde{\kappa} \varrho(y,v);
		\end{eqnarray*}
		which both together reveal that  \eqref{eqIoffeEpigraph} holds true. 
		
	    By Proposition~\ref{propIoffeRegA},  for each $(x,w,y)\in \ball_X(\bar{x}, a) \times \ball_Y(\bar{w},a \ell+c) \times \ball_Y(\bar{z}+\bar{q},b)$ we have
		$\dist((x,w),\mathcal{E}^{-1}(y))\leq \hat{\kappa} \dist \left(y, \mathcal{E}(x,w)\cap\ball_Y[\bar{z} + \bar{q},r]\right)$.
        Fix any  $x\in \ball_X(\bar{x},a)$ and any $y\in\ball_Y(\bar{z}+\bar{q},b)$. Pick an arbitrary $v\in (F+G)(x)\cap \ball_Y[\bar{z}+\bar{q},r]$ (if there is any).	
        As $\varrho(\bar{q}, \bar{w}) \leq b$, {\rm (iii)} yields a $w\in G(x)\cap \ball_Y(\bar{w},a \ell + c)$ such that $v\in F(x)+w$. Then 
 	    \begin{eqnarray*}
			\varrho(v-w,\bar{z})&=& 		\varrho(v-w+\bar{q},\bar{z}+\bar{q})\leq 	\varrho(v-w+\bar{q},v)+	\varrho(v,\bar{z}+\bar{q}) \\
			&  \leq & \varrho(w, \bar{w})+	\varrho(\bar{w}, \bar{q}) + \varrho(v,\bar{z}+\bar{q}) <  a \ell +c + b + r < \delta.
		\end{eqnarray*}
		So $v \in \mathcal{E}(x,w)\cap\ball_Y[\bar{z} + \bar{q},r]$. 
        As  $\dist(x,(F+G)^{-1}(y)) \leq \dist((x,w),\mathcal{E}^{-1}(y))$, we get that 
		$\dist(x,(F+G)^{-1}(y)) \leq \hat{\kappa}\varrho(y,v)$.
		Since $v \in (F+G)(x) \cap \ball_Y[\bar{z} + \bar{q},r]$ is arbitrary, the proof is finished. 
	\end{proof}

An immediate  global consequence is \cite[Theorem~5I.2]{DontchevBook}:  	
\begin{corollary} \label{corEpigraph} Let $(X,d)$ be a  complete metric space and $(Y,\varrho)$ be a linear complete metric space with a translation-invariant metric. Consider constants  $\kappa > 0$ and $\ell > 0$ such that $\kappa \ell < 1$ and set-valued mappings $F$, $G:X\tto Y$ having closed graphs. 
	If $F$ is regular on $X\times Y$ with the constant $\kappa$  and $G$ is Hausdorff-Lipschitz on $X$ with the constant $\ell$,
	then  $F+G$ is regular on $X\times Y$ with any  constant $\hat{\kappa}>\kappa/(1-\kappa\ell)$. 
\end{corollary}
\begin{proof}
	Fix any $\hat{\kappa}>\kappa/(1-\kappa\ell)$. Pick any $(\bar{x}, \bar{y}) \in \gph (F+G)$.    Find points  $\bar{z}\in F(\bar{x})$ and  $\bar{w}\in G(\bar{x})$ such that $\bar{z} + \bar{w} = \bar{y}$.  Given an arbitrary $a > 0$, the assumptions of Theorem~\ref{thmSUSVP_S1}, with $b:= a$, $c := \infty$, $r:=a$, and $\bar{q}:=\bar{w}$, hold. 
 \end{proof}

It is well-known that the third parameter $r$ can be neglected in the conclusion of Theorem~\ref{thmSUSVP_S1} in exchange for a possibly smaller neighborhood of the reference point, cf. \cite[Proposition 5H.1]{DontchevBook}.  

\begin{remark}\label{remConstrainedRegToReg} 
\rm Given a point $(\bar{x},\bar{y})\in X\times Y$ and a  mapping $H:X\tto Y$, with $\bar{y}\in H(\bar{x})$, assume that  there are positive constants   $\delta$, $r$, and $\kappa$  such that $H$ is $\ball_Y(\bar{y},r)-$restrictedly regular on   $\ball_{{X}}(\bar{x},\delta)\times \ball_Y(\bar{y},\delta) =: \Omega$ with the constant $\kappa$. 
Let  $\beta := \min\{\delta, \kappa r/(2 + 2\kappa)\}$.     	Then $H$ is regular on $\ball_X(\bar{x}, \beta) \times \ball_Y(\bar{y}, \beta) = : W$  with the constant $\kappa$.
   	Indeed,  fix any $(x,y) \in W$.  Pick an arbitrary  $v\in H(x)$ (if there is any). Assume that $ \varrho(y,v)\geq r/2$. As  $\beta \leq \delta$, we have  $(\bar{x},y)\in \Omega$. Therefore 
   	\begin{eqnarray*}
   		\dist\big(x,H^{-1}(y)\big) & \leq & 	\dist\big(\bar{x},H^{-1}(y)\big)+d(x,\bar{x})\leq \kappa\varrho(y,\bar{y})+d(x,\bar{x}) <  \beta (1 + \kappa) \\
   		&  \leq &  \kappa r/2  \leq  \kappa\varrho(y,v).
   	\end{eqnarray*}
   	If $\varrho(y,v)<r/2$, then $ \varrho(\bar{y},v) \leq \varrho(\bar{y},y)+\varrho(y,v)<  \beta + r/2 < r$; and, as  $(x,y)\in \Omega$, we get the same inequality immediately.
 \end{remark}

Second, we derive an improved version of \cite[Theorem 3.2]{ACN15}. Note that the original slope-based proof is rather complicated, e.g., the EVP is used twice.

  \begin{corollary} \label{corACN} 
  Let $(X,d)$ be a  metric space and $(Y, \varrho)$ be a linear  metric space with a~translation-invariant metric.  Consider positive constants $\alpha$, $\beta$, $\kappa$,  and $\ell$ such that $\kappa \ell < 1$, a point $(\bar{x}, \bar{z}, \bar{w})\in X\times Y\times Y$, and set-valued mappings $F$, $G:X\rightrightarrows Y$, with $\bar{z}\in F(\bar{x})$ and $\bar{w} \in G(\bar{x})$,  such that $\gph F\cap\big(\ball_X[\bar{x},\alpha]\times \ball_Y[\bar{z},   \beta]\big)$ and $\gph G\cap\big(\ball_X[\bar{x},\alpha]\times \ball_Y[\bar{w},\beta]\big)$ are complete.  If $F$ is regular on $\ball_X(\bar{x},\alpha)\times \ball_Y(\bar{y},\beta)$ with the constant $\kappa$ and $G$ is Hausdorff-Lipschitz on $\ball_X(\bar{x}, \alpha)$ with the constant $\ell$ and such that  $\diam\,G(\bar{x}) < \beta$, then for each $\hat{\kappa}>\kappa/(1-\kappa\ell)$  and each $\delta > 0$ satisfying 
\begin{eqnarray}
	\label{eqConstantsInequlityA}
3\delta( 1 +  \hat{\kappa}) < \alpha \quad \mbox{and} \quad   \delta ( 1 + \hat{\kappa})(3\ell + 2/\hat{\kappa}) + \delta + \diam\,G(\bar{x}) < \beta
\end{eqnarray}
the mapping $F+G$ is regular on  $\ball_X(\bar{x},  \delta) \times \ball_Y(\bar{z} + \bar{w},\delta)$ with the constant $\hat{\kappa}$. 
\end{corollary}

\begin{proof} Let $\hat{\kappa}$ and  $\delta$ be as in the premise.  Let $r := 2\delta(1+\hat{\kappa}) /\hat{\kappa}$ and $\hat{\alpha} := \delta + \hat{\kappa} (\delta + r) = 3\delta (1 +  \hat{\kappa})$. We shall apply Theorem \ref{thmSUSVP_S1},  with  $a:= \delta$, $b:=\delta$, $c := \diam\,G(\bar{x})$, and $\bar{q} := \bar{w}$. As  $\hat{\alpha}  < \alpha$ and $c + \hat{\alpha} \ell + b + r< \beta$, completeness assumptions  as well as (i) and (ii) therein hold true.
	To show {\rm (iii)}, fix any $x\in\ball_X(\bar{x}, \delta )$ and any $v\in (F+G)(x)$. There is a $w\in G(x)$ such that $v\in F(x)+w$. Since $G$ is Hausdorff-Lipschitz on $\ball_X(\bar{x}, \alpha)$,  there is a $\bar{w}'\in G(\bar{x})$ such that $\varrho(w,\bar{w}')\leq\ell \, d(x,\bar{x})$. 
Hence 
\begin{equation} \label{eqWBWrho}
		\varrho(w,\bar{w})\leq \varrho(w,\bar{w}')+\varrho(\bar{w}',\bar{w})\leq \ell \, d(x,\bar{x})+\diam\,G(\bar{x});
	\end{equation}
	so  ${\rm (iii)}$ holds. Theorem \ref{thmSUSVP_S1} and Remark~\ref{remConstrainedRegToReg}, with $\kappa := \hat{\kappa}$,  yield the conclusion because
	 $\min\{ \delta, \hat{\kappa} r /(2+2\hat{\kappa}) \} = \delta$. 
\end{proof}

We get a generalization of \cite[Theorem 5G.3]{DontchevBook} proved as \cite[Theorem 2.3]{cf15}: 

\begin{corollary} \label{cor5G3} 
	Let $(X,d)$ be a  complete metric space and $(Y,\varrho)$ be a linear complete metric space with a translation-invariant metric. Consider positive constants $\alpha$, $\beta$,  $\kappa$, and $\ell$ such that $\kappa\ell<1$, a point $(\bar{x}, \bar{z}) \in X \times Y$, and  a set-valued mapping $F:X\tto Y$, with $\bz\in F(\bx)$, such that $\gph F \cap (\ball_X[\bar x,\alpha] \times\ball_Y[\bz,\beta])$ is 
		closed. Assume that $F$ is $\ball_Y[\bz,\beta]-$restrictedly  regular on $\ball_X(\bx,\alpha)\times\ball_Y(\bz,\beta)$ with the constant $\kappa$.
	Then for each $\hat{\kappa}>\kappa/(1-\kappa\ell)$, each $a > 0$,  each  $b > 0$, and  each $g:X\to Y$, which is Lipschitz continuous on $\ball_X[\bx, a + 2\hat{\kappa} b]$ with the constant $\ell$,  such that  
	\begin{equation}
		\label{eqInequalities2}
				a + 2 \hat{\kappa} b \leq \alpha,\quad  2b +\ell (a + 2 \hat{\kappa} b)\leq \beta, \quad\text{and}\quad 
				\varrho(0,g(\bx)) \leq b,
	\end{equation}			 
	the mapping $F+g$ has the following property:  for each $y$, $\hat{y}\in \ball_Y(\bz,b)$ and each $\hat{x} \in (F+g)^{-1}(\hat{y})\cap \ball_{X}(\bx,a)$ there is an $x \in (F+g)^{-1}(y)\cap  \ball_{X}(\bx,a+ 2 \hat{\kappa} b )$ such that $d(x, \hat{x}) \leq \hat{\kappa} \varrho(y, \hat{y})$.
\end{corollary}
\begin{proof}
	Let $a$, $b$, $\hat{\kappa}$, and $g$ be as in the premise.  Pick  a $\tilde{\kappa} \in (\kappa/(1-\kappa\ell), \hat{\kappa})$. Put $c:=0$, $r := b$, $\bar{w} := g(\bar{x})$, and $\bar{q} := 0$. Define  $\hat{\alpha}:= a + 2\tilde{\kappa}b$, $\hat{\beta}:= \hat{\alpha}\ell$, and $\delta := \hat{\beta} + 2 b$. Then
	$$
		a+ \tilde{\kappa} (b + r) = \hat{\alpha} < \alpha, 
		\quad 
		\hat{\beta} + b + r =  \ell \hat{\alpha} + 2b = \delta < \beta,
		\quad \mbox{and} \quad \bar{q} \in \ball_{Y}[\bar{w},b]. 
	$$
	Applying Theorem~\ref{thmSUSVP_S1}, with $(\alpha, \beta, \hat{\kappa})$ replaced by $(\hat{\alpha}, \hat{\beta}, \tilde{\kappa})$, we get that 
	$F+g$  is $\ball_Y[\bz,b]-$restrictedly  regular on $\ball_X(\bx,a)\times\ball_Y(\bz,b)$ with the constant $\tilde{\kappa}$.
	This implies the desired property because $\tilde{\kappa} < \hat{\kappa}$. 
\end{proof}

Instead of the condition on the diameter  one can assume the so-called local sum-stability introduced by M. Durea and R. Strugariu \cite{DS2012}.

   \begin{definition} \label{defSumStab} Let $(X,d)$ and $(Y,\varrho)$ be metric spaces. Consider a $(\bar{x},\bar{z}, \bar{w})\in X\times Y\times Y$ and  set-valued mappings  $F$, $G: X\rightrightarrows Y$, with ${\bar{z}\in F(\bar{x})}$ and  $\bar{w}\in G(\bar{x})$. The pair $(F,G)$ is said to be \emph{sum-stable around} $(\bar{x},\bar{z},\bar{w})$ if
   	for every $\beta >0$ there is an $\alpha >0$ such that, for every $x\in \ball(\bar{x},\alpha )$ and every $v\in(F+G)(x)\cap \ball_Y(\bar{z}+\bar{w},\alpha),$
   	there exist points  $z\in F(x)\cap \ball_Y(\bar{z},\beta)$ and $w\in G(x)\cap \ball_Y(\bar{w}, \beta)$ such that $v=z + w$.
   \end{definition}
 
Using the above notion, we  obtain \cite[Theorem 8]{HNT2014}:
\begin{corollary} \label{corSUMStableA}
Let $(X,d)$ be a  metric space and $(Y, \varrho)$ be a linear metric space with a~translation-invariant metric.  Consider constants $\kappa >0$ and $\ell >0$ such that $\kappa \ell < 1$,  a point $(\bar{x}, \bar{z}, \bar{w})\in X\times Y\times Y$, and set-valued mappings $F$, $G:X\rightrightarrows Y$, with $\bar{z}\in F(\bar{x})$ and $\bar{w} \in G(\bar{x})$,  such that the sets $\gph F$ and $\gph G$ are complete around $(\bx,\bz)$ and $(\bx,\bar{w})$, respectively.
	If $F$ is regular around  $(\bar{x},\bar{z})$ with the constant $\kappa$, $G$ has Aubin property around $(\bar{x},\bar{w})$ with the constant $\ell$, and the pair $(F,G)$ is sum-stable around $(\bar{x},\bar{z}, \bar{w})$,  then  $F+G$ is regular around  $(\bar{x},\bar{z}+\bar{w})$ with any constant $\hat{\kappa} > \kappa/(1-\kappa \ell)$.
\end{corollary}
\begin{proof}	Let $\hat{\kappa} > \kappa/(1-\kappa \ell)$ be arbitrary.	 Find an  $R > 0$ such that the sets $\gph F\cap\left(\ball_{X}[\bx,R]\times \ball_Y[\bz, R]\right)$ and $\gph G\cap\left(\ball_{X}[\bx,R]\times \ball_Y[\bar{w},R]\right)$ are complete, that $F$ is regular on $\ball_{X}(\bx, R)\times\ball_Y(\bz,R)$ with the constant $\kappa$, and that $G$ has Aubin property on  $\ball_X(\bx,R) \times \ball_Y(\bar{w},R)$ with the constant $\ell$. Pick a $c > 0$ such that 
$$ 
   3c (1+ \hat{\kappa}) < R
   \quad \mbox{and} \quad 
    2c +  3c \ell  (1+ \hat{\kappa})  + 2c (1+ \hat{\kappa})/\hat{\kappa} < R.
$$ 
As  $(F,G)$ is sum-stable around $(\bar{x},\bar{z}, \bar{w})$, there is a $\delta \in (0, c)$ such that for each $x\in\ball_{X}(\bx,\delta)$ and each $v\in (F+G)(x)\cap\ball_Y[\bz+\bar{w}, \delta  + 2\delta(1+ \hat{\kappa})/\hat{\kappa}]$ there is a~$w\in G(x)\cap\ball_Y(\bar{w}, c)$ with  $v\in F(x)+w$. Theorem \ref{thmSUSVP_S1},  with $(a,b,r,\bar{q}) := (\delta,\delta, 2\delta(1+ \hat{\kappa})/\hat{\kappa}, \bar{w})$, and Remark~\ref{remConstrainedRegToReg}, with $\kappa:=\hat{\kappa}$, imply that $F+G$ is regular on $\ball_{X}(\bx,\delta)\times \ball_Y(\bz+\bar{w},\delta)$ with the constant $\hat{\kappa}$ because $\min\{ \delta, \hat{\kappa} r /(2+2\hat{\kappa}) \} = \delta$.
\end{proof}

\subsection{Mappings with a complete domain and a closed graph} \label{sscPScdcg}
In this subsection, we investigate the stability of the second type of regularity under additive perturbations. We start with single-valued perturbations.

\begin{theorem} \label{thmPSonW}
 Let $(X,d)$ be a complete metric space and $(Y,\varrho)$ be a linear metric space with a translation-invariant metric. Consider constants $\kappa > 0$ and $\ell > 0$ such that $\kappa \ell < 1$, a non-empty set $\Omega \subset X\times Y$, a function $\gamma : X \to [0, \infty]$ which is both positive and Lipschitz continuous on $\Omega_X$ with the constant $1$, a set-valued mapping $F:X\tto Y$ with a closed graph, and a single-valued mapping $g:X\to Y$ which is continuous on $X$. Assume that the set   
	\begin{equation} \label{eqPSonW}
	W:=\big\lbrace (x,y)\in X\times Y :  \ \ball_X(x, \gamma(x))\times \ball_Y(y-g(x), \ell \gamma(x))\subset \Omega \big\rbrace.
	\end{equation}
	is such that $\dom F \cap W_X$ is non-empty. 
	If $F$ is $\gamma$-regular on $\Omega$ with the constant $\kappa$  and $g$ is Lipschitz continuous on $\Omega_X$ with the constant $\ell$, then $F + g$ is $\gamma$-regular on $W$ with any constant $\hat{\kappa}>\kappa/(1-\kappa\ell)$.
\end{theorem}
\begin{proof}  Let $\hat{\kappa} > \kappa/(1-\kappa \ell)$ be arbitrary.
	Pick  a constant $c \in (\ell, 1/\kappa)$ such that $ 1/(c - \ell) = \hat{\kappa}$. Find a $\lambda > 1$  such that 
	$ \lambda c < 1/\kappa$ and let $\tilde{\kappa}:= \hat{\kappa}/\lambda$. Then  $\kappa < 1/(\lambda c)  <  \tilde{\kappa} <  \hat{\kappa}$. 
	
	Let $\omega$ be a (compatible) metric on $X \times Y$ defined by \eqref{defOmega} with $\kappa := \tilde{\kappa}$ and  let $H= F + g$. Then $H$ has a closed graph. We check the assumptions of Corollary~\ref{corIoffeGraph1}, with $F:=H$ and $\kappa := \tilde{\kappa}$. To do so, note that  $ \hat{\kappa}/\tilde{\kappa}= \lambda$ and fix  any $u\in X$ and any $y \in W_Y $ such that  $0<\lambda\dist((u,y),\gph H)< \gamma(x) -d(u,{x})$ for an ${x}\in W_{\cdot,y}$.
	We have to find a $\hat{u} \in X$ such that 
	\begin{equation} \label{eqBetteru}
	 \lambda  \dist \big( (\hat{u},y), \gph H \big) <  \lambda  \dist \big( (u,y), \gph H \big) 	 -   d(\hat{u},u).
	\end{equation}
	Pick a pair $(\hat{x},v)\in \gph H$ such that 
	\begin{eqnarray}
		\label{eqSemiInequality2}
		0<d(u,\hat{x})+\tilde{\kappa} \varrho(y,v)< \lambda  \dist\big((u,y),\gph H\big). 
	\end{eqnarray}	
	If $v = y$, then $\hat{u}:=\hat{x}$ satisfies \eqref{eqBetteru} because from \eqref{eqSemiInequality2} we get that
	\begin{eqnarray*}
		\lambda  \dist \big( (\hat{u},y), \gph H \big) & = & 0=d(u, \hat{x}) + \tilde{\kappa} \varrho(v,y) - d(\hat{u},u) \\
		&  <  & \lambda  \dist\big((u,y),\gph H \big)-d(\hat{u},u).
	\end{eqnarray*}
	Assume that $y\neq v$. Since $x \in W_X \subset \Omega_X$, \eqref{eqSemiInequality2} yields that 
	\begin{equation}
		d(\hat{x},x) \leq  	 d(\hat{x},u)+	 d(u,x)< \lambda  \dist\big((u,y),\gph H \big) + d(u,x)
		 <   \gamma(x). \label{eqSemiInequality} 
	\end{equation}
	So $\hat{x}\in \ball_{{X}}(x,\gamma(x)) \subset \Omega_X$.    The Lipschitz continuity of $g$ and \eqref{eqSemiInequality} imply that   
	\begin{eqnarray*}
		\varrho(y-g(\hat{x}),y-g({x})) &=& \varrho(g(\hat{x}),g(x))\leq \ell \, d(\hat{x},x) <   \ell \gamma(x).
	\end{eqnarray*}
	Thus  $y-g(\hat{x})\in \ball_{{Y}}(y-g({x}), \ell  \gamma(x))$. Remembering that $(x,y) \in W$, we conclude that $(\hat{x},y-g(\hat{x}))\in \Omega$. 
    Since $v - g(\hat{x}) \in F(\hat{x})$, \eqref{eqSemiInequality2} implies that 
	\begin{eqnarray*}
		\tilde{\kappa} \dist(y-g(\hat{x}),F(\hat{x})) & \leq &  \tilde{\kappa} \varrho(y,v) < \lambda  \dist\big((u,y),\gph H\big) -  d(u,\hat{x}) \\
		& < &   \gamma(x) -\big( d(u,x) + d(u,\hat{x}) \big) \leq \gamma(x) -d(\hat{x},x) \leq  \gamma(\hat{x}).
	\end{eqnarray*}
	As $\kappa < 1/(c \lambda) < \tilde{\kappa}$, $\gamma$-regularity of $F$ on $\Omega$ yields  a~$\hat{u}\in F^{-1}(y-g(\hat{x}))$ such that 
	\begin{eqnarray*}
		d(\hat{x}, \hat{u})&<&1/(c \lambda)\varrho(y,v) \quad \quad \big( < \gamma(x) -d(\hat{x},x)\big).
	\end{eqnarray*}
	Then $d(x,\hat{u}) \leq d(x,\hat{x}) + d(\hat{x}, \hat{u}) < \gamma(x)$. 	Hence $\hat{u} \in \ball_{{X}}(x, \gamma(x)) \subset \Omega_X$.	 Let $\hat{v}:=y-g(\hat{x})+g(\hat{u})\in H(\hat{u})$. Then $\varrho (y, \hat{v}) = \varrho(g(\hat{x}),g(\hat{u}) \leq \ell \,  d(\hat{x},\hat{u})$; therefore, remembering that $ (c-  \ell) \lambda \tilde{\kappa} = 1$, we get  that
	\begin{eqnarray*}
		\lambda \, \dist \big( (\hat{u},y), \gph H \big) & \leq &  \lambda \, (d(\hat{u}, \hat{u}) +  \tilde{\kappa} \, \varrho (y, \hat{v} )) = \lambda \tilde{\kappa} \,  \varrho (y, \hat{v})\leq  \lambda \tilde{\kappa} \ell \, d(\hat{x},\hat{u}) \\
		& = &  \tilde{\kappa} \, \varrho(y,v) -   \tilde{\kappa} \, \varrho(y,v) +  \lambda \tilde{\kappa} \ell \, d(\hat{x},\hat{u}) \\
		& < &  \tilde{\kappa} \,  \varrho(y,v) - (c-  \ell) \lambda \tilde{\kappa} \,  d(\hat{x},\hat{u})\\
		& = &   d(u,\hat{x}) +  \tilde{\kappa} \, \varrho(y,v) -  \big(  d(u,\hat{x})+d(\hat{x},\hat{u})\big)\\
		& \leq &   d(u,\hat{x}) +  \tilde{\kappa} \, \varrho(y,v) -  d(u,\hat{u}) \\
		& \overset{\eqref{eqSemiInequality2}}{<} & \lambda \, \dist \big( (u,y), \gph H \big) 	 -   d(\hat{u},u).
	\end{eqnarray*}	
	Thus  \eqref{eqBetteru} holds. By Corollary~\ref{corIoffeGraph1},  
	$H$ is $\gamma$-regular on $W$ with the constant $\hat{\kappa}$.
\end{proof}

We immediately obtain the stability of Milyutin regularity. 
   
\begin{theorem} \label{thmPSonWofMilR}
 Let $(X,d)$ be a complete metric space and $(Y,\varrho)$ be a linear metric space with a translation-invariant metric. Consider positive constants $r$, $\varepsilon$, $\kappa$,  and $\ell$ such that $\kappa \ell < 1$, non-empty open sets $U \subset X$ and $V \subset Y$, a set-valued mapping $F:X\tto Y$ with a closed graph, and a continuous single-valued mapping $g:X\to Y$ such that the set 
 \begin{equation} \label{eqPSonWofMilR}
	W_{\varepsilon}=\big\lbrace (x,y)\in X\times Y :  \  x \in U \mbox{ and } \ball_Y(y-g(x), \varepsilon \ell \dist(x, X\setminus U)) \subset V \big\rbrace
	\end{equation}
	is non-empty. 
	If $F$ is Milyutin regular on $U \times V$ with the constants $r$ and $\kappa$ and  $g$ is Lipschitz continuous on $U$ with the constant $\ell$, then for each $\hat{\kappa}>\kappa/(1-\kappa\ell)$ the mapping $F + g$ is Milyutin regular on $W_{\varepsilon}$ with the constants  $\delta := \min\{\kappa r, \varepsilon,1\}/\hat{\kappa}$ and $\hat{\kappa}$.
\end{theorem}

\begin{proof}
 Let $\hat{\kappa}>\kappa/(1-\kappa\ell)$ be arbitrary. Let $\delta := \min\{\kappa r, \varepsilon,1\}/\hat{\kappa}$ and $\gamma(x):= \hat{\kappa} \delta \dist(x, X\setminus U)$, $x \in X$. As $F$ is Milyutin regular on $U \times V$ with the constants $r$ and  $\kappa$ and $\hat{\kappa} \delta \leq \kappa r $, it is $\gamma$-regular on $\Omega := U \times V$ with the constant $\kappa$. Then $\Omega_X = U$.  Let $W$  be defined by \eqref{eqPSonW}.
As $\hat{\kappa} \delta \leq 1$, we have $\ball_X(x, \gamma(x)) \subset \ball_X(x, \dist(x, X\setminus U)) \subset U$ for each $x \in U$. Hence 
$$
  W= \big\lbrace (x,y)\in X\times Y :  \  x \in U \mbox{ and } \ball_Y\big(y-g(x), \hat{\kappa} \delta \ell \dist(x,X\setminus  U)\big)\subset V \big\rbrace.
$$
Theorem~\ref{thmPSonW} implies that $F + g$ is $\gamma$-regular on $W$ with the constant $\hat{\kappa}$. Note that  $W_{\varepsilon} \subset W$ because $\hat{\kappa} \delta \leq \varepsilon$ and that $\dist(x, X \setminus (W_{\varepsilon})_X) \leq \dist(x, X \setminus U)$ for each $x \in X$ because $(W_{\varepsilon})_X \subset U$. Therefore $F + g$ is Milyutin regular on $W_{\varepsilon}$ with the constants $\delta$ and $\hat{\kappa}$.  
\end{proof}

The above statement is a quantitative version of \cite[Theorem 4.2]{NDH2022}; and if  $r := 1 /\kappa$ and $\varepsilon = 1$,  it reduces to \cite[Theorem 4.3]{NDH2022}.  Two \emph{separate} proofs of these statements in \cite{NDH2022} are based on two \emph{different} slope sufficient conditions (with two independent proofs again). Despite this, both the proofs are slightly longer, request the use of sequences, and calculating several limits, etc. On the contrary, we use very definitions of the properties of $F$, $g$, and $\varrho$; the triangle inequality; and the criterion. 

Now, we prove an analogue of Theorem~\ref{thmSUSVP_S1} having several seemingly different reincarnations again.   

	\begin{theorem} \label{thmPSonWSetValuedSemilocalB} Let $(X,d)$ be a complete metric space and $(Y,\varrho)$ be a linear metric space with a translation-invariant metric. Consider positive constants  $\delta$, $\kappa$, $\hat{\kappa}$, and $\ell$ such that $\kappa \ell < 1$ and $\hat{\kappa}>\kappa/(1-\kappa\ell)$, numbers $\mu$, $r \in (0, \infty]$, a point $(\bar{x}, \bar{z}, \bar{w}) \in X \times Y \times Y$, and set-valued mappings $F$, $G:X\tto Y$, with $\bar{z}\in F(\bar{x})$ and $\bar{w} \in G(\bar{x})$, such that the sum $H:=F+G$ has a closed graph. Assume that $\alpha := \delta + \hat{\kappa} r$ is such that 
		\begin{itemize}
			 \item[$(i)$] $F$ is regular with the constant $\kappa$ on  $ \{(x,y) \in \ball_{{X}}(\bar{x}, \alpha) \times \ball_{{Y}}(\bar{z}, \mu + \delta) : \ \dist\big(y,F(x)\cap \ball_Y(\bar{z}, \mu + \delta+r) \big)  < r\}$  ;
			 \item[$(ii)$] $G$ has Aubin property on $\ball_{{X}}(\bar{x}, \alpha) \times \ball_Y(\bar{w}, \mu)$ with the constant $\ell$; 
			 \item[$(iii)$] for each $x \in \ball_{{X}}(\bar{x}, \alpha)$ and each $v \in H(x) \cap \ball_{{Y}}(\bar{z} + \bar{w}, \delta + r)$ there is a $w \in G(x) \cap \ball_Y(\bar{w}, \mu)$ such that $v \in F(x) + w$. 
		 \end{itemize}
		  Then $H$ is regular on $ \{(x,y) \in \ball_X(\bar{x},\delta) \times \ball_Y(\bar{z} + \bar{w},\delta) : \ \dist(y,H(x)) < r \}$ with the constant $\hat{\kappa}$. \end{theorem}

\begin{proof}  
	Fix any $\hat{r} \in (0,r)$. Let $ W:=\ball_X(\bar{x},\delta) \times \ball_Y(\bar{z} + \bar{w},\delta)$ and  $\gamma(x):= \hat{\kappa} \hat{r}$, $x \in W_X = \ball_X(\bar{x},\delta)$.
	
	Pick  a  $c \in (\ell, 1/\kappa)$ such that $ 1/(c - \ell) = \hat{\kappa}$. Find a $\lambda > 1$  such that $\lambda \hat{r} < r$ and
	$ \lambda c < 1/\kappa$. Let $\tilde{\kappa}:= \hat{\kappa}/\lambda$. Then  $\kappa < 1/(\lambda c)  <  \tilde{\kappa} <  \hat{\kappa}$. 
	Let $\omega$ be a (compatible) metric on $X \times Y$ defined by \eqref{defOmega} with $\kappa := \tilde{\kappa}$.  
	
	We check the assumptions of Corollary~\ref{corIoffeGraph1}, with $F:=H$ and  $\kappa := \tilde{\kappa}$.   To do so, note that  $\hat{\kappa} / \tilde{\kappa}= \lambda$ and fix  any $u\in X$ and any $y \in W_Y $ such that  $0<\lambda\dist((u,y),\gph H)< \hat{\kappa} \hat{r} -d(u,{x})$ for an ${x}\in W_{\cdot,y}$.
	We have to find a $\hat{u} \in X$ such that 
	\begin{equation} \label{eqBetteruSetValuedSemilocalTOM}
		\lambda  \dist \big( (\hat{u},y), \gph H \big) <  \lambda  \dist \big( (u,y), \gph H \big) 	 -   d(\hat{u},u).
	\end{equation}
	Pick a pair $(\hat{x},v)\in \gph H$ such that 
	\begin{eqnarray}
		\label{eqSemiInequality2SetValuedSemilocalTOM}
		0<d(u,\hat{x})+\tilde{\kappa} \varrho(y,v)< \lambda  \dist\big((u,y),\gph H\big). 
	\end{eqnarray}	
	If $v = y$, then $\hat{u}:=\hat{x}$ satisfies \eqref{eqBetteruSetValuedSemilocalTOM} because from \eqref{eqSemiInequality2SetValuedSemilocalTOM} we get that
	\begin{eqnarray*}
		\lambda  \dist \big( (\hat{u},y), \gph H \big) & = & 0=d(u, \hat{x}) + \tilde{\kappa} \varrho(v,y) - d(\hat{u},u) \\
		&  <  & \lambda  \dist\big((u,y),\gph H\big)-d(\hat{u},u).
	\end{eqnarray*}
	Assume that $y\neq v$. 
	Inequality \eqref{eqSemiInequality2SetValuedSemilocalTOM} yields that 
	\begin{eqnarray}
		\label{eqSemiInequality3SetValuedSemilocalTOM}
		d(\hat{x},x) \leq  	 d(\hat{x},u)+	 d(u,x)< \lambda  \dist\big((u,y),\gph H\big)+ d(u,x)  <  \hat{\kappa} \hat{r} 
	\end{eqnarray}
and
\begin{eqnarray} \label{eqSetValuedSemilocalTOM}
   \tilde{\kappa} \varrho(y,v) 	& < &   \lambda  \dist\big((u,y),\gph H \big) - d(u,\hat{x}) 
		 <      \hat{\kappa}\hat{r} - d(x,u) - d(u,\hat{x}) \notag \\ 
		 & \leq & \hat{\kappa}\hat{r}  - d(x,\hat{x}) \leq \hat{\kappa}\hat{r}.
\end{eqnarray}
    	As $(x,y) \in W = \ball_X(\bar{x},\delta) \times \ball_Y(\bar{z} + \bar{w},\delta)$ and $\hat{r} < \lambda \hat{r}  < r$, \eqref{eqSemiInequality3SetValuedSemilocalTOM} and \eqref{eqSetValuedSemilocalTOM} divided by $\tilde{\kappa}$ imply that $\hat{x}\in \ball_{{X}}(\bar{x}, \alpha )$ and $v \in \ball_Y(\bar{z}+\bar{w},r+\delta)$, respectively. By {\rm (iii)}, there is a $\hat{w} \in G(\hat{x})\cap\ball_Y(\bar{w},\mu)$ such that $v\in F(\hat{x})+ \hat{w}$.
	Then 
	\begin{eqnarray*}
		\varrho(y-\hat{w},\bar{z}) &=& \varrho(y + \bar{w} - \hat{w} ,\bar{z} + \bar{w}) \leq  \varrho(y + \bar{w} - \hat{w} ,y) + \varrho(y ,\bar{z} + \bar{w}) \\
		& = & \varrho(\bar{w}, \hat{w}) + \varrho(y ,\bar{z} + \bar{w}) < \mu + \delta.
	\end{eqnarray*} 
	Similarly, we get that
	\begin{eqnarray*}
		\varrho(v-\hat{w},\bar{z}) &=& \varrho(v + \bar{w} - \hat{w} ,\bar{z} + \bar{w}) \leq  \varrho(v + \bar{w} - \hat{w} ,v) + \varrho(v ,\bar{z} + \bar{w}) \\
		& = & \varrho(\bar{w}, \hat{w}) + \varrho(v ,\bar{z} + \bar{w}) < \mu + r+ \delta.
	\end{eqnarray*} 
	Since $v - \hat{w} \in F(\hat{x})$,  \eqref{eqSetValuedSemilocalTOM} implies that 
	\begin{eqnarray*}
		\tilde{\kappa} \dist\big(y-\hat{w},F(\hat{x})\cap \ball_Y(\bar{z}, \mu + \delta + r) \big)  & \leq &    \tilde{\kappa} \varrho(y-\hat{w}, v - \hat{w}) =   \tilde{\kappa} \varrho(y,v)\\ 
		& < & \hat{\kappa}\hat{r}  - d(x,\hat{x}) \leq \hat{\kappa}\hat{r}.
	\end{eqnarray*} 
	So $\dist\big(y-\hat{w},F(\hat{x})\cap \ball_Y(\bar{z}, \mu + r+ \delta) \big)  < \lambda \hat{r} < r$. Therefore  $(\hat{x},y-\hat{w})$ lies in the set in (i).  
	As $\kappa < 1/(c \lambda) < \tilde{\kappa}$, employing {\rm (i)} and \eqref{eqSetValuedSemilocalTOM}, we find  a $\hat{u}\in F^{-1}(y-\hat{w})$ such that 
	\begin{eqnarray*}
		d(\hat{x}, \hat{u})&<&1/(c \lambda)\varrho(y,v)  \quad \quad \big( < \hat{\kappa} r -d(\hat{x},x) \big).
	\end{eqnarray*}
	Then $d(x,\hat{u}) \leq d(x,\hat{x}) + d(\hat{x}, \hat{u}) < \hat{\kappa} r$. 	Hence $\hat{u} \in \ball_{{X}}(\bar{x}, \alpha)$. 
	Using {\rm (ii)}, we find a $\tilde{w}\in G(\hat{u})$ such that $\varrho(\hat{w}, \tilde{w})\leq \ell \, d(\hat{x},\hat{u})$. Let $\hat{v}:=y-\hat{w}+\tilde{w} \in H(\hat{u})$. Then $\varrho (y, \hat{v}) = \varrho(\tilde{w},\hat{w}) \leq \ell \,  d(\hat{x},\hat{u})$; therefore, remembering that $ (c-  \ell) \lambda \tilde{\kappa} = 1$, we get  
	\begin{eqnarray*}
		\lambda \, \dist \big( (\hat{u},y), \gph H \big) & \leq &  \lambda \, (d(\hat{u}, \hat{u}) +  \tilde{\kappa} \, \varrho (y, \hat{v} )) = \lambda \tilde{\kappa} \,  \varrho (y, \hat{v})\leq  \lambda \tilde{\kappa} \ell \, d(\hat{x},\hat{u}) \\
		& = &  \tilde{\kappa} \, \varrho(y,v) -   \tilde{\kappa} \, \varrho(y,v) +  \lambda \tilde{\kappa} \ell \, d(\hat{x},\hat{u}) \\
		& < &  \tilde{\kappa} \,  \varrho(y,v) - (c-  \ell) \lambda \tilde{\kappa} \,  d(\hat{x},\hat{u})\\
		& = &   d(u,\hat{x}) +  \tilde{\kappa} \, \varrho(y,v) -  \big(  d(u,\hat{x})+d(\hat{x},\hat{u})\big)\\
		& \leq &   d(u,\hat{x}) +  \tilde{\kappa} \, \varrho(y,v) -  d(u,\hat{u}) \\
		& \overset{\eqref{eqSemiInequality2SetValuedSemilocalTOM}}{<} & \lambda \, \dist \big( (u,y), \gph H \big) 	 -   d(\hat{u},u).
	\end{eqnarray*}	
	Thus \eqref{eqBetteruSetValuedSemilocalTOM} holds.  By Corollary~\ref{corIoffeGraph1}, 
	$H$ is $\gamma$-regular on $W$ with the constant $\hat{\kappa}$. Since $\hat{r} \in (0, r)$ is arbitrary, we obtain the desired conclusion. 
\end{proof}

First, we get a statement parallel to Corollary~\ref{corEpigraph} covering \cite[Corollary 3.4]{HT08}:
  
\begin{corollary} Let $(X,d)$ be a complete metric space and $(Y,\varrho)$ be a linear metric space with a translation-invariant metric. Consider  constants $\kappa > 0$ and $\ell>0$ such that $\kappa \ell < 1$, an $r \in (0, \infty]$,  and set-valued mappings $F$, $G:X\tto Y$ such that the sum $H:=F+G$ has a closed graph. 
	If $F$ is regular on $\{(x,y) \in X\times Y : \ \dist(y,F(x)) < r\}$ with the constant $\kappa$  and $G$ is Hausdorff-Lipschitz on $X$ with the constant $\ell$,
	then  $H$ is regular on $\{(x,y) \in X \times Y: \ \dist(y,H(x)) < r \}$ with any  constant $\hat{\kappa}>\kappa/(1-\kappa\ell)$. 
\end{corollary}
\begin{proof} Let  $\hat{\kappa}>\kappa/(1-\kappa\ell)$ be arbitrary. Pick any $(\bar{x}, \bar{y}) \in \gph H$.    Find points  $\bar{z}\in F(\bar{x})$ and  $\bar{w}\in G(\bar{x})$ such that $\bar{z} + \bar{w} = \bar{y}$. Given an arbitrary $\delta>0$, all the assumptions of Theorem \ref{thmPSonWSetValuedSemilocalB}, with $\mu := \infty$, are satisfied. 
\end{proof}

It is well-known that the second parameter $r$ can be neglected in the conclusion of Theorem~\ref{thmPSonWSetValuedSemilocalB} in exchange for a possibly smaller neighborhood of the reference point, cf. \cite[Exercise 2.12]{IoffeBook}.  
   \begin{remark} \label{remMetricRegularity}  \rm Given a point $(\bar{x},\bar{y})\in X\times Y$ and a  mapping $H:X\tto Y$, with $\bar{y}\in H(\bar{x})$, assume that  there are positive constants  $r$, $\delta$, and $\kappa$  such that $H$ is regular on 
   $	 \lbrace (x,y)\in\ball_{{X}}(\bar{x},\delta)\times \ball_Y(\bar{y},\delta) : \ \dist(y,H(x))<r\rbrace =: \Omega$ with the constant $\kappa$. Let $\beta := \min \{ \delta, \kappa r/(1 + \kappa)\}$.
   	Then $H$ is regular on $\ball_X(\bar{x}, \beta) \times \ball_Y(\bar{y}, \beta) = : W$  with the constant $\kappa$.
   	Indeed,  fix any $(x,y) \in W$.  Pick an arbitrary  $v\in H(x)$ (if there is any). Assume that $ \varrho(y,v)\geq r$. As $\dist(y,H(\bar{x})) \leq \varrho (y, \bar{y}) < \beta < r$ and $\beta \leq \delta$, we have  $(\bar{x},y)\in \Omega$. Therefore 
   	\begin{eqnarray*}
   		\dist\big(x,H^{-1}(y)\big) & \leq & 	\dist\big(\bar{x},H^{-1}(y)\big)+d(x,\bar{x})\leq \kappa\varrho(y,\bar{y})+d(x,\bar{x}) <  \beta (1 + \kappa) \\
   		&  \leq &  \kappa r  \leq  \kappa\varrho(y,v).
   	\end{eqnarray*}
   	If $\varrho(y,v)<r$, then $(x,y)\in \Omega$; and we get the same inequality immediately. 
   \end{remark}
   
 Second, we derive an analogue of  Corollary~\ref{corACN} in the setting $(S_2)$. 

\begin{corollary} \label{corACN_B} Let $(X,d)$ be a complete metric space and $(Y, \varrho)$ be a linear  metric space with a~translation-invariant metric.  Consider positive constants $\alpha$, $\beta$, $\kappa$,  and $\ell$ such that $\kappa \ell < 1$, a point $(\bar{x}, \bar{z}, \bar{w})\in X\times Y\times Y$, and set-valued mappings $F$, $G:X\rightrightarrows Y$, with $\bar{z}\in F(\bar{x})$ and $\bar{w} \in G(\bar{x})$,  such that the sum $H:=F+G$ has a closed graph. If $F$ is regular on $\ball_X(\bar{x},\alpha)\times \ball_Y(\bar{y},\beta)$ with the constant $\kappa$ and $G$ is Hausdorff-Lipschitz on $\ball_X(\bar{x}, \alpha)$ with the constant $\ell$ and such that  $\diam\,G(\bar{x}) < \beta$, then for each $\hat{\kappa}>\kappa/(1-\kappa\ell)$  and each $\delta > 0$ satisfying 
\begin{eqnarray}
	\label{eqConstantsInequlity}
\delta( 2 +  \hat{\kappa}) < \alpha \quad \mbox{and} \quad  \ell \delta ( 2 + \hat{\kappa})+ \delta + \diam\,G(\bar{x}) < \beta
\end{eqnarray}
the mapping $H$ is regular on  $\ball_X(\bar{x},  \delta) \times \ball_Y(\bar{z} + \bar{w},\delta)$ with the constant $\hat{\kappa}$. 
\end{corollary}

\begin{proof} Let $\hat{\kappa}$ and  $\delta$ be as in the premise.  Let $r := \delta(1+\hat{\kappa}) /\hat{\kappa}$, $\hat{\alpha} := \delta + \hat{\kappa} r = \delta (2 +  \hat{\kappa})$,  and $\mu:=\ell \hat{\alpha} + \diam\,G(\bar{x})$. We shall apply Theorem \ref{thmPSonWSetValuedSemilocalB} with $\alpha := \hat{\alpha}$. As  $\hat{\alpha}  < \alpha$, assumption (ii) therein holds trivially; while (i) is satisfied because we also have  $ \mu  + \delta < \beta$.
	To show {\rm (iii)}, fix any $x\in\ball_X(\bar{x},\hat{\alpha})$ and any $v\in H(x)$. There is a $w\in G(x)$ such that $v\in F(x)+w$. Since $G$ is Hausdorff-Lipschitz on $\ball_X(\bar{x}, \alpha)$,  there is a $\bar{w}'\in G(\bar{x})$ such that $\varrho(w,\bar{w}')\leq\ell \, d(x,\bar{x})$. 
By \eqref{eqWBWrho},  ${\rm (iii)}$ holds. Theorem \ref{thmPSonWSetValuedSemilocalB} and Remark~\ref{remMetricRegularity} yield the conclusion because  $\min\{ \delta, \hat{\kappa} r /(1+\hat{\kappa}) \} = \delta$.
\end{proof}
If $G(\bar{x}) = \{\bar{w}\}$ we cover  \cite[Corollary 4.5]{HT08}. Note that Corollary 3.4 and Corollary 4.5 in \cite{HT08} are consequences of two general approximation results with two \emph{separate} proofs based on two \emph{different} slope sufficient conditions (with two independent proofs again); and request also $F$ to have a closed graph  and $Y$ to be a normed space. The following elementary example,  see also \cite[Example 5I.1]{DontchevBook} for an excessively elaborated one, shows a well-known fact that the sum of two set-valued mappings, one of which is locally regular while  the latter is locally Hausdorff-Lipschitz, does not need to be regular around the reference point; in particular, \cite[Corollary 4.8]{HT08} is wrong.    

\begin{example} \label{exSumDiameter} \rm 
 Let $X$, $Y$, $d$, $\varrho$, $F$, and $(\bar{x},\bar{z})$ be as in Example~\ref{exFregIandII}. Let $\bar{w}:=0$. Consider a mapping $\mathbb{R}\ni x \longmapsto G(x) := \lbrace 0,1 \rbrace\subset \mathbb{R}$. Then $F$, $G$, and $H : = F + G$ have closed graphs, $\bar{z} \in F(\bar{x})$,  $\bar{w} \in G(\bar{x})$,  $G$ is Hausdorff-Lipschitz on $X$ with any constant $\ell > 0$, and $H(x)=\{x,-1,1+x,0\}$ for each $x \in \mathbb{R}$.
However $H$ is not regular around $(\bar{x}, \bar{z} +\bar{w}) = (0,0)$. Indeed, let $\hat{\kappa} > 0$, and open sets $U \ni 0$  and $V \ni 0$ be \emph{arbitrary}. Pick a  $t > 0$ such that $y := - t \in V \cap (-1/2,0)$ and $x:=\hat{\kappa}   t\in U$. Then $x > 0$ and $H^{-1}(y) = \{-t, -t-1\}$. Hence 
$$
  \dist(x, H^{-1}(y)) = \hat{\kappa}   t + t  > \hat{\kappa}   t  =   \hat{\kappa} |y - 0| = \hat{\kappa} \dist(y,H(x)).  
$$
Remember that given a $\beta > 1 = \diam G(\bar{x})$, there are no $\alpha > 0$ and $\kappa > 0$ such that $F$ is regular on $(-\alpha, \alpha) \times (-\beta, \beta)$ with the constant $\kappa$.  Hence the condition on the diameter in Corollary~\ref{corACN} and Corollary~\ref{corACN_B} cannot be improved. 
\end{example}

Finally, we obtain an analogue of Corollary~\ref{corSUMStableA}, that is, of \cite[Theorem 8]{HNT2014} in the setting $(S_2)$.

\begin{corollary} \label{corSUMStableB}
	Let $(X,d)$ be a complete metric space and $(Y, \varrho)$ be a linear  metric space with a~translation-invariant metric.  Consider constants $\kappa >0$ and $\ell >0$ such that $\kappa \ell < 1$,  a point $(\bar{x}, \bar{z}, \bar{w})\in X\times Y\times Y$, and set-valued mappings $F$, $G:X\rightrightarrows Y$, with $\bar{z}\in F(\bar{x})$ and $\bar{w} \in G(\bar{x})$,  such that the sum $H:=F+G$ has a closed graph.
	If $F$ is regular around  $(\bar{x},\bar{z})$ with the constant $\kappa$, $G$ has Aubin property around $(\bar{x},\bar{w})$ with the constant $\ell$, and the pair $(F,G)$ is sum-stable around $(\bar{x},\bar{z}, \bar{w})$,  then $H$ is regular around  $(\bar{x},\bar{z}+\bar{w})$ with any constant $\hat{\kappa} > \kappa/(1-\kappa \ell)$.
\end{corollary}
\begin{proof}	Let $\hat{\kappa} > \kappa/(1-\kappa \ell)$ be arbitrary.	Find a $\beta >0$  such that $F$ is regular on $\ball_X(\bar{x},\beta) \times \ball_{{Y}}(\bar{z},2\beta)$ with the constant $\kappa$ and that $G$ has Aubin property on $\ball_X(\bar{x},\beta) \times \ball_Y(\bar{w}, \beta)$ with the constant $\ell$.  
	Since the pair $(F, G)$ is sum-stable around $(\bar{x},\bar{z},\bar{w})$, there is an $\hat{\alpha} \in (0, \beta)$ such that for each $x\in\ball_{{X}}(\bar{x},\hat{\alpha})$ and each $v\in H(x)\cap \ball_{{Y}}(\bar{z}+\bar{w}, \hat{\alpha})$, there is a $w\in G(x)\cap \ball_Y(\bar{w},\beta)$ such that $v\in F(x)+w$. Pick a $\delta \in (0, \beta)$ such that 
	$$ 
	  \delta (2 + \hat{\kappa} ) < \hat{\alpha}
	  \quad \mbox{and} \quad 
	  \delta (1 + 2\hat{\kappa}) <  \hat{\kappa}\hat{\alpha}.
	$$  
	Putting $\mu := \beta$, $r := \delta (1 + \hat{\kappa})/\hat{\kappa}$, and $\alpha := \delta + \hat{\kappa} r$, we have
	$$
	 \alpha < \hat{\alpha} < \beta, \quad 
	 \delta + r < \hat{\alpha},
	 \quad \mbox{and} \quad 
	 \mu + \delta < 2 \beta.
	$$
	Theorem \ref{thmPSonWSetValuedSemilocalB} and Remark~\ref{remMetricRegularity} yield that $H$ is regular on $\ball_{{X}}(\bar{x},\delta) \times \ball_{{Y}}(\bar{z}+\bar{w}, \delta)$ with the constant $\hat{\kappa}$.
\end{proof}

\section{Concluding remarks} \label{scCR}
We briefly discuss possible generalizations and topics for further research:  
\paragraph*{\bf 1.} In \cite{NDH2022}, non-linear analogues of the properties in Definition~\ref{defRegBtypes} are considered. Both Proposition~\ref{propIoffeRegA} and Proposition~\ref{propIoffeGraph} can be modified accordingly. Since these modifications are \emph{easy and immediate}, we leave them as an exercise for an interested reader (if any);
\paragraph*{\bf 2.} Applying Corollary~\ref{corIoffeSemiRegA} one immediately gets (variants of)  \cite[Theorem 4.4]{CFK2018}, which, in particular, implies that the sum of two set-valued mappings, one of which is locally regular while the other has Aubin property, is \emph{semiregular} at the reference point, cf. Corollary~\ref{corSUMStableA}; while Corollary~\ref{corIoffeGraphSemireg} yields a similar statement in the setting $(S_2)$; 
\paragraph*{\bf 3.} Using a technique similar to the proof of Theorem~\ref{thmSUSVP_S1} and Theorem~\ref{thmPSonWSetValuedSemilocalB} one can investigate
regularity (as well as directional versions of it)  for  compositions as in \cite[Theorem 16]{CDPS} under the assumption of the so-called \emph{composition stability} \cite[Definition 14]{CDPS}, see also \cite{DS14} and \cite{DHNS14}. Note that the slope-based proofs in \cite{DHNS14} are more complicated than the ones of more general results in \cite{CDPS} based on the Ioffe-type criterion;    
\paragraph*{\bf 4.} One can also investigate a stronger version of the regularity than (II) considered in \cite{Ioffe10}, where \eqref{defMRB} is replaced by 
  \begin{equation} \label{defMRC}
		\dist \big(x,F^{-1}(y)\cap W_{\cdot,y} \big)\leq  \kappa \dist \big(y, F(x)\big)
		\quad \mbox{for each} \quad (x,y) \in W.
	\end{equation}
	It covers, in particular, directional regularity introduced in \cite{AAI07}. Theorem 13 in \cite{Ioffe10} guarantees the perturbation stability of this property under single-valued additive perturbations in the setting $(S_1)$. The proof of it is based on the Ioffe-type criterion. In \cite[Theorem 8]{HT15}, the same is established in the setting $(S_2)$ via slope based conditions and lower semi-continuous envelopes. It seems that our approach works also in this case. However, in view of the next point, we prefer to leave a detailed discussion for another note (if any);
\paragraph*{\bf 5.} One can consider the set $W$ in the last displayed inequality as the graph of another mapping $G: X \rightrightarrows Y$. We arrive at the (equivalent) property: {\it for each  $x \in \dom G$ we have that}   
  \begin{equation} \label{defCoiP}
  	\dist \big(x,F^{-1}(y)\cap G^{-1}(y) \big)\leq  \kappa \inf\{ \varrho(y,v): \ (y,v)\in G(x) \times F(x) \}.
  \end{equation}
  In particular, if this holds for an  $x \in  \dom G \cap \dom F$, then the set of the \emph{coincidence points} $\Xi:=\{\xi \in X: \ F(\xi)\cap G(\xi) \neq \emptyset\}$ of $F$ and $G$ is non-empty and we have an estimate to the distance from $x$ to $\Xi$. In \cite{AAGDO09}, results on coincidence points of two (set-valued) mappings can be found, which can be proved either by the iterative procedure or by using the EVP. Let us emphasize that this approach, in contrast to Theorem~\ref{thmSUSVP_S1},  requests that one of the mappings $F$ and $G$ has a complete graph while the other has a \emph{closed} graph  only, see also \cite[Theorem 3]{DF2011} and compare it with Corollary~\ref{corEpigraph}.  Roughly speaking,  the criterion for the property in \eqref{defMRC} is obtained by applying the EVP to the function  $ \gph F \ni (u,v) \longmapsto \varphi (u,v) := \varrho(y, v) + \dist((u,y), W)$, with a given $y \in W_Y$, provided that $F$ has a complete graph \cite{Ioffe10}. If $W := \gph G$, then the roles of $F$ and $G$ can be interchanged which does not effect \eqref{defCoiP}.   However, any serious discussion on the relation between these two concepts takes several pages, see also the umbrella theorem \cite[Theorem 6.1]{Ioffe11} corresponding to the setting $(S_1)$.

\end{document}